\documentclass{amsart}
\usepackage{amsmath,amssymb}
\usepackage{enumerate}
\usepackage{xspace}
\usepackage{boxedminipage}
\usepackage{listings}
\usepackage{tabularx}
\usepackage{subfigure}
\usepackage{algpseudocode}
\usepackage{tristan}
\usepackage{times}

\numberwithin{equation}{section}
\setboolean{showtodo}{false}
\setboolean{shownotes}{false}
\setboolean{showchanges}{false}
\setboolean{usemathrsfs}{true}
\renewcommand{\highlight}{}

\author{Tristan Pryer}
\address{
   Tristan Pryer,
   \thanks{
     Department of Mathematics and Statistics,
     Whiteknights,
     PO Box 220,
     Reading RG6 6AX,
     UK
  {\tt{T.Pryer@reading.ac.uk}}.
}}

\title[dG methods for the $p$-biharmonic equation]{Discontinuous Galerkin methods for the $p$--biharmonic equation from a discrete variational perspective}
\date{\today}
\pdfformat{false}
\begin{document}

\maketitle
\begin{abstract}
  We study discontinuous Galerkin approximations of the
  $p$--biharmonic equation \highlight{for $p\in (1,\infty)$
    \margnote{Ref A, comment (6)}} from a variational perspective. We
  propose a discrete variational formulation of the problem based on a
  appropriate definition of a finite element Hessian and study
  convergence of the method (without rates) using a weak lower
  semicontinuity argument. We also present numerical experiments aimed at
  testing the robustness of the method.
\end{abstract}

\section{Introduction, problem setup and notation}
\label{sec:intro}

The $p$--biharmonic equation is a fourth order elliptic boundary value
problem, related to, in fact a nonlinear generalisation of, the
biharmonic problem. Such problems typically arise from areas of
elasticity\highlight{, in particular the nonlinear case can be used as a model
for travelling waves in suspension bridges
\cite{LazerMcKenna:1990,GyulovMorosanu:2010}. \margnote{Ref A, comment (22}} It is a fourth order
analog to its second order sibling, the $p$--Laplacian, and, as such, is
useful as a prototypical nonlinear fourth order problem.

The efficient numerical simulation of general fourth order problems
has attracted recent interest. A conforming approach to this class of
problem would require the use of $\cont{1}$ finite elements, the
Argyris element for example \cite[Section 6]{Ciarlet:1978}. From a
practical point of view the approach presents difficulties, in that
the $\cont{1}$ finite elements are difficult to design and complicated
to implement, especially when working in three spatial dimensions.

Discontinuous Galerkin (dG) methods form a class of nonconforming
finite element method. They are extremely popular due to their
successful application to an ever expanding range of problems. A very
accessible unification of these methods together with a detailed
historical overview is presented in
\cite{ArnoldBrezziCockburnMarini:2001}.

If $p=2$ we have the special case that the ($2$--)biharmonic problem
is linear. It has been well studied in the context of dG methods, for
example, the papers \cite{LasisSuli:2003, GeorgoulisHouston:2009}
study the use of $h$--$p$ dG finite elements (where $p$ here means the
local polynomial degree) applied to the ($2$--)biharmonic problem. To
the authors knowledge there is currently no finite element method
posed for the general $p$--biharmonic problem.

In this work we use discrete variational techniques to build a
discontinuous Galerkin (dG) numerical scheme for the $p$--biharmonic
operator \highlight{with $p\in (1,\infty)$ \margnote{Ref A, comment
    (6)}}. We are interested in such a methodology due to the
applications to discrete symmetries, in particular, discrete versions
of Noether's Theorem \cite{Noether:1971}.

A key constituent to the numerical method for the $p$--biharmonic
problem (and second order variational problems in general) is an
appropriate definition of the Hessian of a piecewise smooth
function. To formulate the general dG scheme for this problem from a
variational perspective one must construct an appropriate notion of a
Hessian of a piecewise smooth function. The \emph{finite element
  Hessian} was first coined by \cite{Aguilera:2008} for use in the
characterisation of discrete convex functions. Later in
\cite{LakkisPryer:2011a} it was used in a method for nonvariational
problems where the strong form of the PDE was approximated and
put to use in the context of fully nonlinear problems in
\cite{LakkisPryer:2013}. 

Convergence of the method we propose is proved using the framework set
out in \cite{Di-PietroErn:2010} where some extremely useful discrete
functional analysis results are given. Here, the authors use the
framework to prove convergence for a dG approximation to the steady
state incompressible Navier--Stokes equations. A related but
independent work containing similar results is given in
\cite{BuffaOrtner:2009} where the authors study dG approximations to
generic first order variational minimisation problems.

The rest of the paper is set out as follows: The rest of this section
introduces necessary notation and the model problem we consider. In
Section \ref{sec:cont} we give some properties of the continuous
$p$--biharmonic problem. In Section \ref{sec:discretisation} we give the
methodology for discretisation of the model problem. In
Section \ref{sec:convergence} we detail solvability and the convergence of
the discrete problem. Finally, in Section \ref{sec:numerics} we study the
discrete problem computationally and summarise numerical experiments.

Let $\W \subset\reals^d$ be a bounded domain with boundary $\partial
\W$. We begin by introducing the Sobolev spaces
\cite{Ciarlet:1978,Evans:1998}
\begin{gather}
  \leb{p}(\W)
  =
  \ensemble{\phi}
           {\int_\W \norm{\phi}^p < \infty}   \text{ for } p\in[1,\infty) 
             \text{ and }
             \leb{\infty}(\W)
             =
             \ensemble{\phi}
                      {\esssup_\W \norm{\phi} < \infty},
                      \\
                      \sob{l}{p}(\W) 
                      = 
                      \ensemble{\phi\in\leb{p}(\W)}
                               {\D^{\vec\alpha}\phi\in\leb{p}(\W), \text{ for } \norm{\geovec\alpha}\leq l}
                               \text{ and }
    \sobh{l}(\W)
    := 
    \sob{l}{2}(\W),
\end{gather}
which are equipped with the following norms and semi-norms:
\begin{gather}
  \Norm{v}_{\leb{p}(\W)}^p
  :=
       {\int_\W \norm{v}^p}
       \\
       \Norm{v}_{l,p}^p
       := 
       \Norm{v}_{\sob{l}{p}(\W)}^p
       = 
       \sum_{\norm{\vec \alpha}\leq k}\Norm{\D^{\vec \alpha} v}_{\leb{p}(\W)}^p
       \\
       \norm{v}_{l,p}^p
       :=
       \norm{v}_{\sob{l}{p}(\W)}^p
       =
       \sum_{\norm{\vec \alpha} = k}\Norm{\D^{\vec \alpha} v}_{\leb{p}(\W)}^p
       \\
       \Norm{v}_l^2 := \Norm{v}_{\sobh{l}(\W)}^2 = \Norm{v}_{\sob{l}{2}(\W)}^2,
\end{gather}
where $\vec\alpha = \{ \alpha_1,\dots,\alpha_d\}$ is a
multi-index, $\norm{\vec\alpha} = \sum_{i=1}^d\alpha_i$ and
derivatives $\D^{\vec\alpha}$ are understood in a weak sense. We pay
particular attention to the cases $l = 1,2$ and
\begin{gather}
  \sobz{2}{p}(\W) := \ensemble{\phi\in\sob{2}{p}(\W)}{\phi = \Transpose{\qp{\nabla \phi}}\geovec n = 0}.
\end{gather}

In this paper we use the convention that the derivative $\D u$ of a
function $u:\W\to\reals$ is a row vector, while the gradient of $u$,
$\nabla u$ is the derivatives transpose, i.e., $\nabla u =
\Transpose{\left(\D u\right)}$. We will make use of the slight abuse
of notation, following a common practice, whereby the Hessian of $u$
is denoted as $\Hess u$ (instead of the correct $\nabla \D u$) and
is represented by a $d\times d$ matrix.

Let $L = L\qp{\geovec x, u, \nabla u, \Hess u}$ be the \emph{Lagrangian}. We
will let
\begin{equation}
  \label{eq:action-functional}
  \dfunkmapsto[]
  {\cJ[\ \cdot \ ; p]}
  \phi
  {\sobz{2}{p}(\W)}
  {\cJ[\phi; p] := \int_\W L(\geovec x, \phi, \nabla \phi, \Hess \phi) \d \geovec x}
  {\reals}
\end{equation}
be known as the \emph{action functional}. For the $p$--biharmonic
problem the action functional is given explicitly as
\begin{equation}
  \cJ[u;p]
  :=
  \int_\W L(\geovec x, u, \nabla u, \Hess u)
  =
  \int_\W \frac{1}{p}\norm{\Delta u}^p - f u,
\end{equation}
where $\Delta u := \trace\qp{\Hess u}$ is the Laplacian and
$f\in\leb{q}(\W)$ is a known source function. We then look to find a
minimiser over the space $\sobz{2}{p}(\W)$, that is, to find
$u\in\sobz{2}{p}(\W)$ such that
\begin{equation}
  \cJ[u;p] = \min_{v\in\sobz{2}{p}(\W)} \cJ[v;p].
\end{equation}

If we assume temporarily that we have access to a smooth minimiser,
\ie $u\in\cont{4}(\W)$, then, given that the Lagrangian is of second
order, we have that the Euler--Lagrange equations are (in general)
fourth order.

Let $\frob{\geomat X}{\geomat Y} = \trace\qp{\Transpose{\geomat X}
  \geomat Y}$ be the Frobenious inner product between matrices. We
then let
\begin{equation}
  \geomat X = \squmatmd{x}{d}
\end{equation}
then use
\begin{equation}
  {\frac{\partial L}{\partial\qp{\geomat X}}}
  :=
  \squmatmd{{\partial L}/{\partial{x}}}{d}.
\end{equation}
The Euler--Lagrange equations for this problem then take the following
form:
\begin{equation}
  \label{eq:Euler-Lagrange-equations-4th-order}
  \cL[u; p] :=
  \frob{\Hess}{\qp{\frac{\partial L}{\partial\qp{\Hess u}}}}
  +
  \frac{\partial L}{\partial u} = 0.
\end{equation}
These can then be calculated to be
\begin{equation}
  \label{eq:p-biharm}
  \cL[u; p] :=
  \Delta\qp{\norm{\Delta u}^{p-2}\Delta u} - f = 0.
\end{equation}
Note that, for $p=2$, the problem coincides with the biharmonic
problem $\Delta^2 u = f$ which is well studied in the context of dG
methods \cite[e.g.]{Baker:1977, SuliMozolevski:2007, GudiNatarajPani:2008,
  GeorgoulisHouston:2009}.

\section{Properties of the continuous problem}
\label{sec:cont}

To the authors knowledge the numerical method presented here is the
first finite element method presented for the $p$--biharmonic
problem. As such, we will state some simple properties of the problem
which are well known for the problem's second order counterpart \emph{the
$p$--Laplacian} \cite{Ciarlet:1978, BurmanErn:2008}.

\highlight{
\begin{Pro}[equivalence of norms over $\sobz{2}{p}(\W)$ {\cite[Cor 9.10]{Gilbarg:1983}}]
\label{prop:equiv-of-norms}
  Let $\W$ be a bounded domain with Lipschitz boundary then the norms
  $\Norm{\cdot}_{2,p}$ and $\Norm{\Hess \cdot}_{\leb{p}(\W)}$ are
  equivalent over $\sobz{2}{p}(\W)$. \margnote{Ref A, comment (1)}
\end{Pro}}

\begin{Pro}[coercivity of $\cJ$]
  \label{prop:coercive}
  Let $u\in\sobz{2}{p}(\W)$ and $f\in\leb{q}(\W)$, where
  $\tfrac{1}{p}+\tfrac{1}{q}=1$, we have that the action functional
  $\cJ[\ \cdot \ ;p]$ is coercive over $\sobz{2}{p}(\W)$, that is,
  \begin{equation}
    \cJ[u; p] 
    \geq 
    C
    \norm{u}^p_{2,p}
    -\gamma,
  \end{equation}
  for some $C >0 \AND \gamma \geq 0$. Equivalently, let 
  \begin{equation}
    \bi{u}{v}{p} = \int_\W \norm{\Delta u}^{p-2} \Delta u \Delta v
  \end{equation}
  then we have that there exists a constant $C>0$ such that
  \begin{equation}
    \label{eq:coercive-bilinear}
    \bi{v}{v}{p} \geq
    C
    \norm{v}_{2,p}^p \Foreach v\in\sobz{2}{p}(\W).
  \end{equation}
\end{Pro}

\begin{Proof}
  By definition of the $\sobz{2}{p}(\W)$ norm \highlight{and
    Proposition \ref{prop:equiv-of-norms}} we have that
  \begin{equation}
    \cJ[u; p] \geq C(p)\norm{u}_{2,p}^p - fu.
  \end{equation}
  Upon applying H\"older 
  \highlight{and Poincar\'e--Friedrichs inequalities \margnote{Ref A, comment (16)}}
  we see
  \begin{equation}
    \begin{split}
      \cJ[u; p] 
      &\geq
      C(p)
      \norm{u}_{2,p}^p 
      -
      \Norm{f}_{\leb{q}(\W)}\Norm{u}_{\leb{p}(\W)}
      \\
      &\geq
      C(p)\norm{u}_{2,p}^p 
      - C\Norm{f}_{\leb{q}(\W)}.
    \end{split}
  \end{equation}
  The statement (\ref{eq:coercive-bilinear}) is clear 
  \highlight{due to Proposition \ref{prop:equiv-of-norms}}, thus
    concluding the proof.
\end{Proof}

\begin{Pro}[convexity of $L$]
  \label{pro:convex}
  The Lagrangian of the $p$--biharmonic problem is convex with respect
  to its fourth argument.
\end{Pro}
\begin{Proof}
  Using similar arguments to \cite[Section 5.3]{Ciarlet:1978} (also found
  in \cite{BarrettLiu:1994}) the convexity of the functional $J$ is
  a consequence of the convexity of the mapping 
  \begin{equation}
    \cF : \xi\in\reals \to \frac{1}{p}\Norm{\xi}^p.
  \end{equation}
  
\end{Proof}

\begin{Cor}[weak lower semicontinuity]
  \label{cor:semicont}
  The action functional $\cJ$ is weakly lower semicontinuous
  over $\sobz{2}{p}(\W)$. That is, given a sequence of functions $\{
  u_j \}_{j\in\naturals}$ who has a weak limit $u\in\sobz{2}{p}(\W)$, then
  \begin{equation}
    \cJ[u; p] \leq \liminf_{j\to\infty} \cJ[u_j; p].
  \end{equation}
\end{Cor}
\begin{Proof}
  The proof of this is a straightforward extension of \cite[Section 8.2 Thm
    1]{Evans:1998} to second order Lagrangians, noting that
  $\cJ$ is coercive (from Proposition \ref{prop:coercive}) and
  that $L$ is convex with respect to its fourth variable (from
  Proposition \ref{pro:convex}). We omit the full details for brevity.
\end{Proof}

\begin{Cor}[existence and uniquness]
  \label{cor:unique}
  There exists a unique minimiser to the $p$--biharmonic
  equation. Equivalently there is a unique (weak) solution to the (weak)
  Euler--Lagrange equations, find $u\in\sobz{2}{p}(\W)$ such that
  \begin{equation}
    \int_\W \norm{\Delta u}^{p-2} \Delta u \Delta \phi 
    =
    \int_\W f \phi 
    \Foreach \phi\in\sobz{2}{p}(\W).
  \end{equation}
\end{Cor}
\begin{Proof}
  Again, the result can be deduced by extending the arguments in
  \cite[Section 8.2]{Evans:1998} or \cite[Thm 5.3.1]{Ciarlet:1978}, again,
  noting the results of Propositions \ref{prop:coercive} and
  \ref{pro:convex}. The full argument is omitted for brevity.
\end{Proof}

\renewcommand{\vec}[1]{\geovec{#1}}

\section{Discretisation}
\label{sec:discretisation}

Let $\T{}$ be a conforming, shape regular triangulation of $\W$,
namely, $\T{}$ is a finite family of sets such that
\begin{enumerate}
\item $K\in\T{}$ implies $K$ is an open simplex (segment for $d=1$,
  triangle for $d=2$, tetrahedron for $d=3$),
\item for any $K,J\in\T{}$ we have that $\closure K\meet\closure J$ is
  a full subsimplex (i.e., it is either $\emptyset$, a vertex, an
  edge, a face, or the whole of $\closure K$ and $\closure J$) of both
  $\closure K$ and $\closure J$ and
\item $\union{K\in\T{}}\closure K=\closure\W$.
\end{enumerate}
The shape regularity of $\T{}$ is defined as the number
\begin{equation}
  \label{eqn:def:shape-regularity}
  \mu(\T{}) := \inf_{K\in\T{}} \frac{\rho_K}{h_K},
\end{equation}
where $\rho_K$ is the radius of the largest ball contained inside
$K$ and $h_K$ is the diameter of $K$. An indexed family of
triangulations $\setof{\T n}_n$ is called \emph{shape regular} if 
\begin{equation}
  \label{eqn:def:family-shape-regularity}
  \mu:=\inf_n\mu(\T n)>0.
\end{equation}
We use the convention where $\funk h\W\reals$ denotes the \highlight{piecewise constant} \emph{meshsize function} of $\T{}$, i.e.,
\begin{equation}
  h(\vec{x}):=\max_{\closure K\ni \vec x}h_K,
\end{equation}
which we shall commonly refer to as $h$.

We let $\E{}$ be the skeleton (set of common interfaces) of the
triangulation $\T{}$ and say $e\in\E$ if $e$ is on the interior of
$\W$ and $e\in\partial\W$ if $e$ lies on the boundary $\partial\W$
\highlight{and set $h_e$ to be the diameter of $e$.}

\highlight{ We also make the assumption that the mesh is sufficiently
  shape regular such that for any $K\in\T{}$ we have the existence of
  a constant such that
  \begin{equation}
    \label{eq:ass-shape-reg}
    \sum_{e\in\partial K} h_e \norm{e} \leq C \norm{K},
  \end{equation}
  where $\norm{e}$ and $\norm{K}$ denote the $d-1$ and $d$ dimensional measure of $e$ and $K$ respectively.
\margnote{Ref B, comment (1)}}

We let $\poly k(\T{})$ denote the space of piecewise polynomials of
degree $k$ over the triangulation $\T{}$,\ie
\begin{equation}
  \poly k (\T{}) = \{ \phi \text{ such that } \phi|_K \in \poly k (K) \}
\end{equation}
 and introduce the \emph{finite element space}
\begin{gather}
  \label{eqn:def:finite-element-space}
  \fes := \dg{k} = \poly k(\T{}) 
\end{gather}
to be the usual space of discontinuous piecewise polynomial
functions.

\begin{Defn}[finite element sequence]
  A finite element sequence $\{v_h, \fes\}$ is a sequence of discrete
  objects, indexed by the mesh parameter $h$, individually represented
  on a particular finite element space, $\fes$, which itself has
  discretisation parameter $h$, that is, we have that $\fes=\fes(h)$.
\end{Defn}

\begin{Defn}[broken Sobolev spaces, trace spaces]
  \label{defn:broken-sobolev-space}
  We introduce the broken Sobolev space
  \begin{equation}
    \sob{l}{p}(\T{})
    :=
    \ensemble{\phi}
             {\phi|_K\in\sob{l}{p}(K), \text{ for each } K \in \T{}}.
  \end{equation}
  We also make use of functions defined in these broken spaces
  restricted to the skeleton of the triangulation. This requires an
  appropriate trace space
  \begin{equation}
    \Tr{\E} := \prod_{K\in\T{}} \leb{2}(\partial K) 
    \highlight{\supset}
    \prod_{K\in\T{}} \sob{l-\frac{1}{2}}{p}(K)
  \end{equation}
  for $p \geq 2$, $l\geq 1$.
\end{Defn}

\begin{Defn}[jumps, averages and tensor jumps]
  \label{defn:averages-and-jumps}
  We may define average, jump and tensor jump operators over $\Tr{\E}$ for
  arbitrary scalar functions $v\in\Tr{\E}$ and vectors $\vec v\in\Tr{\E}^d$.
  \begin{equation}
    \label{eqn:average}
    \dfunkmapsto[]
	        {\avg{\cdot}}
	        v
	        {\Tr{\E\cup \partial\W}}
	        {
                  \begin{cases}
                    \frac{1}{2}\qp{v|_{K_1} + v|_{K_2}} \text{ over } \E
                    \\
                    v|_{\partial\W} \text{ on } \partial\W
                  \end{cases}
                }
	        {\leb{2}(\E\cup \partial\W)}
  \end{equation}
  \begin{equation}
    \label{eqn:average-vec}
    \dfunkmapsto[]
	        {\avg{\cdot}}
	        {\vec v}
	        {\qb{\Tr{\E\cup \partial\W}}^d}
	        {
                  \begin{cases}
                    \frac{1}{2}\qp{\vec{v}|_{K_1} + \vec{v}|_{K_2}}\text{ over } \E
                    \\
                    \vec v|_{\partial\W} \text{ on } \partial\W
                  \end{cases}
                }
	        {\qb{\leb{2}(\E\cup \partial\W)}^d}
  \end{equation}
  \begin{equation}
    \label{eqn:jump}
    \dfunkmapsto[]
	        {\jump{\cdot}}
	        {v}
	        {{\Tr{\E\cup \partial\W}}}
	        {
                  \begin{cases}
                    {{v}|_{K_1} \geovec n_{K_1} + {v}|_{K_2}} \geovec n_{K_2}\text{ over } \E
                    \\
                    \qp{v \vec n}|_{\partial\W} \text{ on } \partial\W
                  \end{cases}
                }
	        {\qb{\leb{2}(\E\cup \partial\W)}^d}
  \end{equation}
  \begin{equation}
    \label{eqn:jump-vec}
    \dfunkmapsto[]
	        {\jump{\cdot}}
	        {\vec v}
	        {\qb{\Tr{\E\cup \partial\W}}^d}
	        {
                  \begin{cases}
                    {\Transpose{\qp{\vec{v}|_{K_1}}}\geovec n_{K_1} 
                      +
                      \Transpose{\qp{\vec{v}|_{K_2}}}\geovec n_{K_2}}\text{ over } \E
                    \\
                    \qp{\Transpose{\vec v} \vec n}|_{\partial\W} \text{ on } \partial\W
                  \end{cases}
                }
	        {{\leb{2}(\E\cup \partial\W)}}
  \end{equation}
  \begin{equation}
    \label{eqn:tensor-jump}
    \dfunkmapsto[.]
	        {\tjump{\cdot}}
	        {\vec v}
	        {\qb{\Tr{\E\cup \partial\W}}^d}
	        {
                  \begin{cases}
                    {\vec{v}|_{K_1} }\otimes \geovec n_{K_1} + \vec{v}|_{K_2} \otimes\geovec n_{K_2}
                    \text{ over } \E
                    \\
                    \qp{{\vec v} \otimes \vec n}|_{\partial\W} \text{ on } \partial\W
                  \end{cases}
                }
	        {\qb{\leb{2}(\E\cup \partial\W)}^{d\times d}}
  \end{equation}
\end{Defn}

We will often use the following Proposition which we state in full for
clarity but whose proof is merely using the identities in Definition
\ref{defn:averages-and-jumps}.
\begin{Pro}[elementwise integration]
  \label{Pro:trace-jump-avg}
  For a generic vector valued function $\geovec p$ and scalar valued
  function $\phi$ we have
  \begin{equation}
    \label{eq:jump-avg-eq1}
    \begin{split}
      \sum_{K\in\T{}}
      \int_K \div\qp{\geovec p} \phi \d \geovec x
      =
      \sum_{K\in\T{}}
      \qp{
        -
        \int_K
        \Transpose{\geovec p} \nabla_h \phi \d \geovec x
        +
        \int_{\partial K}
        \phi \Transpose{\geovec p} \geovec n_K \d s
        }.
    \end{split}
  \end{equation}
  In particular, if we have $\geovec p \in \Tr{\E\cup\partial\W}^d$ and $\phi \in
  \Tr{\E\cup\partial\W}$, the following identity holds
  \begin{equation}
    \label{eq:jump-avg-eq2}
    \sum_{K\in\T{}}
    \int_{\partial K}
    \phi \Transpose{{\geovec p}}
    \geovec n_K \d s
    =
    \int_\E 
    \jump{\geovec p} 
    \avg{\phi}
    \d s
    +
    \int_{\E\cup\partial\W}
    \Transpose{\jump{\phi}}
    \avg{\geovec p}
    \d s
    =
    \int_{\E\cup\partial\W}
    \jump{\geovec p \phi}
    \d s.
  \end{equation}
  An equivalent tensor formulation of
  (\ref{eq:jump-avg-eq1})--(\ref{eq:jump-avg-eq2}) is
  \begin{equation}
    \label{eq:jump-avg-eq3}
    \begin{split}
      \sum_{K\in\T{}}
      \int_K \D_h {\geovec p} \phi \d \geovec x
      =
      \sum_{K\in\T{}}
      \qp{
        -
        \int_K
        {\geovec p} \otimes \nabla_h \phi \d \geovec x
        +
        \int_{\partial K}
        \phi {\geovec p} \otimes \geovec n_K \d s
        }.
    \end{split}
  \end{equation}
  In particular the following identity holds
  \begin{equation}
    \label{eq:jump-avg-eq4}
    \sum_{K\in\T{}}
    \int_{\partial K}
    \phi {\geovec p}
    \otimes
    \geovec n_K \d s
    =
    \int_\E 
    \tjump{\geovec p} 
    \avg{\phi}
    \d s
    +
    \int_{\E\cup\partial\W}
    {\jump{\phi}}\otimes
    \avg{\geovec p}
    \d s
    =
    \int_{\E\cup\partial\W}
    \tjump{\geovec p \phi}
    \d s.
  \end{equation}
\end{Pro}

The discrete problem we then propose is to minimise an appropriate
discrete action functional, that is to seek $u_h\in\fes$ such that
\begin{equation}
  \cJ_h[u_h; p] = \inf_{v_h\in\fes}\cJ_h[v_h; p].
\end{equation}

\begin{Rem}[motivation for discrete action functional]
  The choice of discrete action functional is crucial. A naive choice
  would be to take the piecewise gradient \highlight{and Hessian} operators,
  substituting them directly into the Lagrangian, \ie
  \begin{equation}
    \cJ_h[u_h; p] 
    =
    \int_\W L\qp{\geovec x, u_h, \nabla_h u_h, \Hess_h u_h}.
  \end{equation}
  This is, however, an inconsistent notion of the derivative operators
  (as noted in \cite{BuffaOrtner:2009}).

  Since for the biharmonic problem the Lagrangian is only dependant on
  the Hessian of the sought function, we need only construct an
  appropriate consistent notion of discrete Hessian.
\end{Rem}

\begin{The}[dG Hessian]
  \label{the:fully-generalised-fe-hessian}
  \highlight{
    Let $v\in\sobz{2}{p}(\T{})$, $\widehat v : \sobh1(\T{}) \to
    \Tr{\E\cup\partial\W}$ be a linear form and $\widehat{\geovec p} :
    \sobh2(\T{}) \times \sobh1(\T{})^d \to \Tr{\E\cup\partial\W}^d$ a
    bilinear form representing consistent numerical fluxes, \ie
  \begin{equation}
    \widehat v \qp{v} = v\vert_{\E\cup\partial\W} 
    \qquad
    \widehat p(v, \nabla v) = \nabla v\vert_{\E\cup\partial\W},
  \end{equation}
  in the spirit of \cite{ArnoldBrezziCockburnMarini:2001}.  Then
  the we define the dG Hessian, $\H[v] \in \fes^{d\times d}$,
  to be the $\leb{2}$ \emph{Reisz representor} of the distributional Hessian of $v$. This has the general form \margnote{Ref A, comments (7,8)}}
  \begin{equation}
    \begin{split}
      \int_{\W} \H[v] \ \Phi
      &=
      -
      \int_\W \nabla_h v \otimes \nabla_h \Phi
      -
      \int_{\E \cup \partial \W} \jump{\widehat v - v} \otimes \avg{\nabla_h \Phi}
      \\
      &\qquad- 
      \int_{\E} \avg{\widehat v - v} \tjump{\nabla_h \Phi}
      +
      \int_{\E\cup \partial\W} \jump{\Phi}\otimes\avg{\geovec{\widehat p}}
      +
      \int_\E \avg{\Phi} \tjump{\geovec{\widehat p}}
      \\
      &\qquad\qquad \Foreach \Phi\in\fes.
    \end{split}
  \end{equation} 
\end{The}

\begin{Proof}
Note that, in view of Green's Theorem, for smooth functions,
$w\in\cont{2}(\W)\cap\cont{1}(\closure{\W})$, we have
\begin{equation}
  \int_\W \Hess w \phi
  =
  -\int_\W \nabla w\otimes \nabla \phi
  +
  \int_{\partial\W} \nabla w \otimes \geovec n \phi
  \Foreach \phi \in\cont{1}(\W)\cap \cont{0}(\closure{\W}).
\end{equation}

As such for a broken function $v\in\sobz{2}{p}(\T{})$ we introduce an
auxiliary variable $\geovec p = \nabla_h v$ and consider the following
primal form for the representation of the Hessian of said function:
For each $K\in \T{}$
\begin{gather}
  \label{eq:primal-hessian}
  \int_K {\H[v]} \ {\Phi} 
  =
  -
  \int_K \geovec p \otimes \nabla_h \Phi
  +
  \int_{\partial K} \geovec{\widehat{p}} \otimes \geovec n \ \Phi
  \Foreach \Phi \in \fes
  \\
  \label{eq:primal-hessian-2}
  \int_K \geovec p \otimes \geovec q 
  =
  -
  \int_K v \ \D \geovec q 
  +
  \int_{\partial K} \geovec q \otimes \geovec n  \ \widehat{v} 
  \Foreach \geovec q \in \fes^d,
\end{gather}
where $\nabla_h = \Transpose{\qp{\D_h }}$ is the elementwise spatial
gradient.

Noting the identity \eqref{eq:jump-avg-eq4} and taking the sum of
(\ref{eq:primal-hessian}) over $K\in\T{}$ we see
\begin{equation}
  \begin{split}
    \int_\W {\H[v]} \ {\Phi} 
    &=
    \sum_{K\in\T{}}
    \int_K {\H[v]} \ {\Phi} 
    =
    \sum_{K\in\T{}}   
    \qp{-  
    \int_K \geovec p \otimes \nabla_h \Phi
    +
    \int_{\partial K} \geovec{\widehat{p}} \otimes \geovec n \ \Phi}
    \\
    &=
    -  
    \int_\W \geovec p \otimes \nabla_h \Phi
    +
    \int_{\E \cup \partial \W} \jump{\Phi} \otimes \avg{\geovec{\widehat{p}}} 
    +
    \int_\E \avg{\Phi} \tjump{\geovec{\widehat{p}}}
  \end{split}
\end{equation}
Using the same argument for (\ref{eq:primal-hessian-2})
\begin{equation}
  \begin{split}
    \int_\W \geovec p \otimes \geovec q 
    &=
    \sum_{K\in\T{}}
    \int_K \geovec p \otimes \geovec q 
    =
    \sum_{K\in\T{}}
    \qp{
      -
      \int_K v \ \D_h \geovec q 
      +
      \int_{\partial K} \geovec q \otimes \geovec n  \ \widehat{v} 
    }
    \\
    &=
    -
    \int_\W v \ \D_h \geovec q 
    +
    \int_{\E \cup \partial \W} \jump{\widehat{v}} \otimes \avg{\geovec q}
    +
    \int_\E \avg{\widehat{v}} \tjump{\geovec q}
  \end{split}
\end{equation}

Note that, again making use of (\ref{eq:jump-avg-eq4}) we have for
each $\geovec q\in\sobh1(\T{})^d$ and $w\in \sobh1(\T{})$ that
\begin{equation}
  \label{eq:edge-identity-2}
  \int_\W \geovec q \otimes \nabla_h w 
  =
  -
  \int_\W \D_h \geovec q w
  +
  \int_{\E\cup \partial \W} \avg{\geovec q}\otimes \jump{w}
  +
  \int_\E \tjump{\geovec q} \avg{w}.
\end{equation}
Taking $w=v$ in (\ref{eq:edge-identity-2}) and substituting into
(\ref{eq:primal-hessian-2}) we see
\begin{equation}
  \label{eq:primal-hessian-3}
  \int_\W \geovec p\otimes \geovec q
  =
  \int_\W \geovec q \otimes \nabla_h v
  +
  \int_{\E \cup \partial \W} \jump{\widehat v - v}\otimes \avg{\geovec q}
  +
  \int_\E \avg{\widehat v - v} \tjump{\geovec q}.
\end{equation}
Now choosing $\geovec q = \nabla_h \Phi$ and substituting
(\ref{eq:primal-hessian-3}) into (\ref{eq:primal-hessian}) concludes
the proof.
\end{Proof}

\begin{Example}[]
  \label{ex:IP-Hess}
  An example of the possible choices of fluxes are
  \begin{gather}
    \widehat v = \begin{cases}
      \avg{v} \text{ over } \E
      \\
      0 \text{ on } \partial \W
      \end{cases}
    \\
    \widehat{\geovec p}
    =
    \avg{\nabla_h v} \text{ on } \E\cup \partial\W.
  \end{gather}
  The result is an interior penalty (IP) type method
  \cite{DouglasDupont:1976} applied to represent the finite element
  Hessian
  \begin{equation}
    \begin{split}
      \int_{\W} \H[v] \ \Phi
      &=
      -
      \int_\W \nabla_h v \otimes \nabla_h \Phi
      +
      \int_{\E \cup \partial \W} \jump{v} \otimes \avg{\nabla_h \Phi}
      \\
      &\qquad
      +
      \int_{\E\cup \partial\W} \jump{\Phi}\otimes\avg{\nabla_h v}.
      \\
      &=
      \int_\W \Hess_h v  \Phi
      -
      \int_{\E\cup\partial\W} \tjump{\nabla_h v} \avg{\Phi}
      \\
      &\qquad
      +
      \int_{\E \cup \partial \W} \jump{v} \otimes \avg{\nabla_h \Phi} 
      .
    \end{split}
  \end{equation}
  \highlight{
    This will be the form of the dG Hessian which we will take for the rest of this exposition.
  }
\end{Example}

\begin{Defn}[lifting operators]
  \label{def:lifting-operators}
  From the IP-Hessian defined in Example \ref{ex:IP-Hess} we define
  the following lifting operator $l_1, l_2:\fes\to\fes^{d\times d}$ such that
  \begin{gather}
    \label{eq:lift-1}
    \int_\W
    l_1[v_h] \Phi
    =
    \int_{\E\cup\partial\W}
    \jump{v_h}
    \otimes
    \avg{\nabla_h \Phi}
    \\
    \label{eq:lift-3}
    \int_\W
    l_2[v_h] \Phi
    =
    - \int_{\E\cup\partial\W}
    \tjump{\nabla_h u_h} \avg{\Phi}.
  \end{gather}
  As such we may write the IP-Hessian as $\H:\fes\to\fes^{d\times d}$ such that 
  \begin{equation}
    \label{eq:definition-of-hessian}
    \int_\W \H[v_h] \Phi = \int_\W \qp{\Hess_h v_h + l_1[v_h] + l_2[v_h]} \Phi \Foreach\Phi\in\fes,
  \end{equation}
  \highlight{where $\Hess_h $ denotes the piecewise Hessian operator.}
\end{Defn}
\highlight{
\begin{Rem}[relation to the local continuous/discontinuous Galerkin method (LCDG)]
  When $\H[\cdot]$ restricted to acting on functions in $\fes\cap\hoz(\W)$ we have that
  \begin{equation}
    \int_\W \H[v_h] \Phi 
    =
    \int_\W \qp{\Hess v_h
    +
    l_2[v_h]} \Phi \Foreach \Phi\in\fes\cap\hoz(\W).
  \end{equation}
  This definition coincides with the auxilliary variable introduced in
  \cite{HuangHuangHan:2010} for Kirchoff plate problems. In addition it is the auxilliary variable used in 
  \cite{LakkisPryer:2011a, LakkisPryer:2013} for applications to
  second order nonvariational PDEs and fully nonlinear PDEs. \margnote{Ref A, comment (9)}
\end{Rem}
}
\section{Convergence}
\label{sec:convergence}
In this section we use the discrete operators from
Section \ref{sec:discretisation} to build a consistent discrete variational
problem and in addition prove convergence. To that end, we being by
defining the natural dG norm for the problem.

\begin{Defn}[dG norm]
  We define the dG norm for this problem as
  \margnote{Ref A, comment (10)}
  \begin{equation}
  \highlight{
    \enorm{v_h}{p}^p
    := 
    \Norm{\Hess_h v_h}^p_{\leb{p}(\W)}
    +
    h\highlight{_e}^{1-p}\Norm{\jump{\nabla_h v_h}}^p_{\leb{p}(\E\cup\partial\W)}
    +
    h\highlight{_e}^{1-2p}\Norm{\jump{v_h}}^p_{\leb{p}(\E\cup\partial\W)},
  }
  \end{equation}
  \margnote{Ref A, comment (1)}
  where \highlight{$\Norm{\cdot}_{\leb{p}(\E\cup\partial\W)}$ is the $d-1$
    dimensional $\leb{p}$ norm over $\E\cup\partial\W$.}\margnote{Ref A, comment (11)}
\end{Defn}

To prove convergence for the $p$-biharmonic equation we modify the
arguments given in \cite{Di-PietroErn:2010} to our problem. To keep
the exposition clear we will, where possible, use the same notation as
in \cite{Di-PietroErn:2010}. 

We state some basic propositions, that is, a trace inequality and
inverse inequality in $\leb{p}(\W)$, the proof of these is readily
available in \cite[e.g.]{Ciarlet:1978}. Henceforth in this section and
throughout the rest of the paper we will use $C$ to denote an
arbitrary positive constant which may depend upon $\mu, p \AND \W$ but
is independent of $h$.

\begin{Pro}[trace inequality]
  \label{pro:trace}
  Let $v_h\in\fes$ be a finite element function then for $p\in
  (1,\infty)$ there exists a constant $C> 0$ such that
  \begin{equation}
    \Norm{v_h}_{\leb{p}({\E\cup\partial\W})}
    \leq
    C h^{-1/p}
    \Norm{v_h}_{\leb{p}(\W)}.
  \end{equation}
\end{Pro}

\begin{Pro}[inverse inequality]
  \label{pro:inverse}
  Let $v_h\in\fes$ be a finite element function then for $p\in
  (1,\infty)$ there exists a constant $C> 0$ such that
  \begin{gather}
    \Norm{\nabla_h v_h}^p_{\leb{p}(\W)}
    \leq
    C h^{-p} 
    \Norm{v_h}^p_{\leb{p}(\W)}
  \end{gather}
\end{Pro}

\begin{Lem}[relating $\enorm{\cdot}{s}$ and $\enorm{\cdot}{t}$ norms]
  \label{eq:equiv-of-dg-norm}
  For two integers $s,t$ such that $1\leq s < t < \infty$ we have that
  there exists a constant $C > 0$ such that 
  \begin{equation}
    \enorm{v_h}{s} 
    \leq
    C
    \enorm{v_h}{t}.
  \end{equation}
\end{Lem}
\begin{Proof}
  The proof follows a similar line to \cite[Lem
  6.1]{Di-PietroErn:2010}. By definition of the $\enorm{\cdot}{s}$
  norm we have that
  \begin{equation}
    \begin{split}
      \enorm{v_h}{s}^s 
      &=
      \int_\W 
      \norm{\Hess_h v_h}^s
      +
      h\highlight{_e}^{1-s} 
      \int_{\E\cup\partial\W}
      \norm{\jump{\nabla_h v_h}}^s
      +
      h\highlight{_e}^{1-2s}
      \int_{\E\cup\partial\W}
      \norm{\jump{v_h}}^s.
    \end{split}
  \end{equation}
  Now let us denote $r = \tfrac{t}{s}$ and $q =
  \tfrac{r}{r-1}$, that is, we have that
  $\tfrac{1}{r}+\tfrac{1}{q} = 1$. Hence we may deduce that
  \margnote{Ref A, comment (12)}
  \begin{equation}
    \begin{split}
      \enorm{v_h}{s}^s 
      &=
      \int_\W \norm{\Hess_h v_h}^s
      +
      \int_{\E\cup\partial\W} h\highlight{_e}^{{1}/{q}} h\highlight{_e}^{{\qp{1-t}}/{r}} \norm{\jump{\nabla_h v_h}}^s
      +
      \int_{\E\cup\partial\W} h\highlight{_e}^{1/q} h\highlight{_e}^{\qp{1-2t}/r} \norm{\jump{v_h}}^s
      \\
      &
      \leq
      \qp{\int_\W 1^{q}}^{1/{q}}
      \qp{\int_\W \norm{\Hess_h v_h}^t}^{1/r} 
      +
      \qp{\highlight{h_e}\int_{\E\cup\partial\W} 1^q }^{1/q}
      \qp{\int_{\E\cup\partial\W} h\highlight{_e}^{1-t} \norm{\jump{\nabla_h v_h}}^t}^{1/r}
      \\
      &\qquad
      +
      \qp{h\highlight{_e}\int_{\E\cup\partial\W} 1^q}^{1/q}
      \qp{\int_{\E\cup\partial\W} h\highlight{_e}^{1-2t} \norm{\jump{v_h}}^t}^{1/r}
      \\
      &
      \leq
      C \enorm{v_h}{t}^s
    \end{split}
  \end{equation}
  where we have used a H\"older inequality together with
  \begin{gather}
    1-s = 1-\tfrac{t}{r} = \tfrac{1}{q} + \tfrac{1- t}{r} \AND
    \\
    1-2s = 1-\tfrac{2t}{r} = \tfrac{1}{q} + \tfrac{1 - 2t}{r},
  \end{gather}
  and the shape regularity of $\T{}$ given in
  (\ref{eq:ass-shape-reg}), concluding the proof.
\end{Proof}

\begin{Defn}[bounded variation]
  Let $\cV[\cdot]$ denote the variation functional defined as
  \begin{equation}
    \cV[u] 
    :=
    \sup
    \ensemble{\int_\W u\div{\geovec \phi}}{\geovec \phi\in[\cont{1}_0(\W)]^d, \Norm{\geovec \phi}_{\leb{\infty}(\W)}\leq 1}.
  \end{equation}
  The space of \emph{bounded variations} denoted \emph{BV} is the
  space of functions with bounded variation functional,
  \begin{equation}
    BV := \ensemble{\phi\in\leb{1}(\W)}{\cV[\phi] < \infty}.
  \end{equation}
  Note that the variation functional defines a norm over $BV$, we set
  \begin{equation}
    \Norm{u}_{BV}
    =
    \cV[u].
  \end{equation}
\end{Defn}

\begin{Pro}[control of the $\leb{\frac{d}{d-1}}(\W)$ norm \cite{EymardGallouetHerbin:2010}]
  \label{pro:control-on-l-1}
  Let $u\in BV$ then we have that there exists a constant $C$ such
  that
  \begin{equation}
    \Norm{u}_{\leb{\frac{d}{d-1}}(\W)} 
    \leq
    C \Norm{u}_{BV}.
  \end{equation}
\end{Pro}

\begin{Pro}[broken Poincar\'e inequality \cite{BuffaOrtner:2009}]
  \label{pro:broken-poincare}
  For $v_h\in\fes$ we have that
  \begin{equation}
    \Norm{v_h}_{\leb{1}(\W)} \leq C \qp{\int_\W \norm{\nabla_h v_h} + \int_{\E\cup\partial\W} \norm{\jump{v_h}}}.
  \end{equation}
\end{Pro}

\begin{Lem}[control on the BV norm]
  \label{lem:control-on-bv}
  We have that for each $v_h\in\fes$ and $p\in[1,\infty)$ that there
    exists a constant $C>0$ such that
  \begin{equation}
    \Norm{v_h}_{BV} \leq C\enorm{v_h}{p}
  \end{equation}
\end{Lem}
\begin{Proof}
  Owing to \cite[Lem 6.2]{Di-PietroErn:2010} we have that
  \begin{equation}
    \label{eq:bv-bound}
    \Norm{v_h}_{BV} \leq \int_\W \norm{\nabla_h v_h} + \int_{\E\cup\partial\W} \norm{\jump{v_h}}.
  \end{equation}
  Applying the broken Poincar\'e inequality given in Proposition
  \ref{pro:broken-poincare} to the first term on the
  (\ref{eq:bv-bound}) gives \margnote{Ref A, comment (1)}
  \begin{equation}
    \begin{split}
      \Norm{v_h}_{BV} 
      &\leq
      C\qp{\int_\W \highlight{\norm{\Hess_h v_h}} + \int_{\E\cup\partial\W} \norm{\jump{\nabla_h v_h}} + \int_{\E\cup\partial\W} \norm{\jump{v_h}}}
      \\
      &\leq
      C\qp{\int_\W \norm{\Hess_h v_h} + \int_{\E\cup\partial\W} \norm{\jump{\nabla_h v_h}} + h\highlight{_e}^{-1}\int_{\E\cup\partial\W} \norm{\jump{v_h}}}
      \\
      &\leq
      C\enorm{v_h}{1}.
    \end{split}
  \end{equation}
  Applying Lemma \ref{eq:equiv-of-dg-norm} concludes the proof.
\end{Proof}

\begin{Lem}[discrete Sobolev embeddings]
  \label{lem:discrete-sob-embed}
  For $v_h\in\fes$ there exists a constant $C>0$ such that
  \begin{equation}
    \Norm{v_h}_{\leb{p}(\W)} 
    \leq 
    C\enorm{v_h}{p}.
  \end{equation}
\end{Lem}
\begin{Proof}
  The proof mimics that of the Gagliardo--Nirenberg--Sobolev
  inequality in \cite[Thm 1, p.263]{Evans:1998}.

  We begin by noting that Proposition \ref{pro:control-on-l-1}
  together with Lemma \ref{lem:control-on-bv} infers the result for
  $p=1$, \ie
  \begin{equation}
    \label{eq:embed-p-1}
    \Norm{v_h}_{\leb{1}(\W)} 
    \leq 
    C\enorm{v_h}{1}.   
  \end{equation}
  Now, we divide the remaining cases into two possibilities, $p\in
  (1,d)$ and $p\in [d,\infty)$. 
    \newline    \newline
    Step 1.
    We begin with $p\in (1,d)$. First note that \highlight{the result of Proposition \ref{pro:control-on-l-1} together with Lemma \ref{lem:control-on-bv} infer that
      \begin{equation}
        \Norm{v_h}_{\leb{\frac{d}{d-1}}(\W)} 
        \leq
        C\enorm{v_h}{1} \Foreach v_h \in\fes.
      \end{equation}
    }
    Now choose $v_h = \norm{w_h}^\gamma$, where $\gamma > 1$ is
    to be chosen, we see
    \begin{equation}
      \label{eq:step-1}
      \qp{\int_\W \norm{w_h}^{\frac{\gamma d}{d-1}}}^{\frac{d-1}{d}} 
      \leq 
      C
      \qp{
        \int_\W \norm{\Hess_h \qp{\norm{w_h}^\gamma}}
      +
      \int_{\E\cup\partial\W} \norm{\jump{\nabla_h \qp{\norm{w_h}^{\gamma}}}}
      +
      \int_{\E\cup\partial\W} h\highlight{_e}^{-1}\norm{\jump{\norm{w_h}^{\gamma}}}}.
    \end{equation}
    We proceed to bound each of these terms individually. \highlight{Firstly note that by the chain rule, we have that \margnote{Ref B, comment (2)}
    \begin{equation}
      \nabla_h\qp{\norm{w_h}^\gamma}
      =
      \gamma \norm{w_h}^{\gamma-1} \nabla_h\qp{\norm{w_h}}
      =
      \gamma \norm{w_h}^{\gamma-2} w_h \nabla_h w_h.
    \end{equation}
    Hence we see that
    \begin{equation}
      \begin{split}
        \Hess_h\qp{\norm{w_h}^\gamma}
        &=
        \D_h \qp{\nabla_h {\norm{w_h}^\gamma}}
        =
        \D_h \qp{\gamma \norm{w_h}^{\gamma-2} w_h \nabla_h w_h}
        \\
        &=
        \gamma \qp{
          \D_h \qp{\norm{w_h}^{\gamma-2}} w_h \nabla_h w_h
          +
          \norm{w_h}^{\gamma-2}\D_h w_h \nabla_h w_h
          +
          \norm{w_h}^{\gamma-2}w_h \Hess_h w_h
        }
        \\
        &=
        \gamma\qp{\gamma-1} \norm{w_h}^{\gamma-2} \nabla_h w_h \otimes \nabla_h w_h
        +
        \gamma \norm{w_h}^{\gamma-2}w_h \Hess_h w_h.
     \end{split}
    \end{equation}
}
\highlight{
Using a triangle
inequality it follows that
\begin{equation}
  \begin{split}
    \int_\W \norm{\Hess_h\qp{\norm{w_h}^\gamma}}
    \leq
    \gamma 
    \int_\W
        \norm{
          \norm{w_h}^{\gamma-1} \Hess_h w_h
        } 
    +
    \gamma\qp{\gamma-1} 
    \int_\W
        \norm{
          \norm{w_h}^{\gamma-2} \nabla_h w_h \otimes \nabla_h w_h
        }
  \end{split}
\end{equation}
}

By a H\"older inequality we have that
    \begin{equation}
      \int_\W \norm{w_h}^{\gamma-1}\norm{\Hess_h w_h}
      \leq
      \qp{\int_\W \norm{w_h}^{q\qp{\gamma - 1}}}^{\frac{1}{q}}
      \qp{\int_\W \norm{\Hess_h w_h}^{p}}^{\frac{1}{p}},
    \end{equation}
    where $q = \frac{p}{p-1}$. 

\highlight{
  In addition we have 
  \begin{equation}
    \int_\W \norm{\norm{w_h}^{\gamma-2} \nabla_h w_h \otimes \nabla_h w_h}
    \leq
    \qp{
      \int_\W 
      \norm{\norm{w_h}^{\gamma-2} \nabla_h w_h}^q
    }^{\frac{1}{q}}
    \qp{
      \int_\W \norm{\nabla_h w_h}^p
    }^{\frac{1}{p}}.
  \end{equation}
Noting that 
\begin{equation}
  \nabla_h\qp{\norm{w_h}^{\gamma-1}} 
  =
  \qp{\gamma-1} \norm{w_h}^{\gamma-3} w_h \nabla_h w_h,
\end{equation}
we see
\begin{equation}
  \begin{split}
  \int_\W \norm{\norm{w_h}^{\gamma-2} \nabla_h w_h \otimes \nabla_h w_h}
  &\leq
  \frac{1}{\gamma-1}
  \qp{
    \int_\W 
    \norm{\nabla_h \qp{\norm{w_h}^{\gamma-1}}}^q
  }^{\frac{1}{q}}
    \qp{
      \int_\W \norm{\nabla_h w_h}^p
    }^{\frac{1}{p}}
    \\
    &\leq
    \frac{C}{\gamma-1}
    \qp{
    \int_\W \norm{w_h}^{q\qp{\gamma-1}}
    }^{\frac{1}{q}}
    \qp{
      \int_\W\norm{\Hess_h w_h}^p
    }^{\frac{1}{p}}
    \end{split}
\end{equation}
by the inverse inequality from Proposition \ref{pro:inverse}.
\margnote{Ref A, comment (2)}}

\highlight{
Hence we have that
\begin{equation}
 \label{eq:bound-on-lap}
 \int_\W \norm{\Hess_h\qp{\norm{w_h}^\gamma}}
  \leq
  C\gamma
    \qp{
    \int_\W \norm{w_h}^{q\qp{\gamma-1}}
    }^{\frac{1}{q}}
    \qp{
      \int_\W\norm{\Hess_h w_h}^p
    }^{\frac{1}{p}}.
\end{equation}
}
    
    Now we must bound the skeletal terms appearing in
    (\ref{eq:step-1}). The jump terms here also act like derivatives
    in that they satisfy a 'chain rule' inequality, using the
    definition of the jump and average operators it holds that
    \begin{equation}
      \label{eq:bound-on-jump-grad}
      \begin{split}
        \int_{\E\cup\partial\W} \norm{\jump{\nabla_h \norm{w_h}^{\gamma}}}
        &\leq
        \int_{\E\cup\partial\W} 2\gamma \avg{\norm{w_h}^{\gamma-1}}\jump{\nabla_h w_h}
        \\
        &\leq
        2\gamma 
        \Norm{h\highlight{_e}^\alpha \avg{\norm{w_h}^{\gamma-1}}}_{\leb{q}({\E\cup\partial\W})}
        \Norm{h\highlight{_e}^{-\alpha} \jump{\nabla_h{w_h}}}_{\leb{p}({\E\cup\partial\W})},
      \end{split}
    \end{equation}
    by a H\"older inequality.

\highlight{
  Focusing our attention to the average term it holds, in view
    of the trace inequality in Proposition \ref{pro:trace}, that
    \begin{equation}
      \begin{split}
        \Norm{h\highlight{_e}^\alpha \avg{\norm{w_h}^{\gamma-1}}}^q_{\leb{q}({\E\cup\partial\W})}
        &\leq 
        C\sum_{K\in\T{}} h\highlight{_e}^{q\alpha - {1}}       
        \Norm{\norm{w_h}^{\gamma-1}}^q_{\leb{q}(K)}
        \\
        &\leq
        C h\highlight{_e}^{q\alpha - {1}}       
        \qp{
          \int_\W \norm{w_h}^{q\qp{\gamma-1}}
        }.
      \end{split}
    \end{equation}
    Upon taking the $q$--th root we see
    \begin{equation}
      \label{eq:bound-on-avg}
      \Norm{h\highlight{_e}^\alpha \avg{\norm{w_h}^{\gamma-1}}}_{\leb{q}({\E\cup\partial\W})}
      \leq
      C h\highlight{_e}^{\alpha - \frac{1}{q}}       
        \qp{
          \int_\W \norm{w_h}^{q\qp{\gamma-1}}
        }^{\frac{1}{q}}.
    \end{equation}
  \margnote{Ref A, comment (3)}}  Choosing $\alpha=\frac{1}{q}$ such that the exponent of $h$
    vanishes and substituting into (\ref{eq:bound-on-jump-grad})
    gives
    \begin{equation}
      \label{eq:final-bound-on-jump-grad}
      \int_{\E\cup\partial\W} \norm{\jump{\nabla_h \norm{w_h}^{\gamma}}}
      \leq
        C 
        \qp{
          \int_\W \norm{w_h}^{q\qp{\gamma-1}}
        }^{\frac{1}{q}}
        \Norm{h\highlight{_e}^{-{\frac{1}{q}}} \jump{\nabla_h{w_h}}}_{\leb{p}({\E\cup\partial\W})}.
    \end{equation}
    
    The final term is dealt with in much the same way. Again, using
    the 'chain rule' type inequality we see that
    \begin{equation}
      \begin{split}
        \int_{\E\cup\partial\W} h\highlight{_e}^{-1} \norm{\jump{\norm{w_h}^{\gamma}}}
        &\leq
        2\gamma \int_{\E\cup\partial\W}
        h\highlight{_e}^{-1}\avg{\norm{w_h}^{\gamma-1}}\norm{\jump{w_h}}
        \\
        &\leq
        2\gamma 
        \Norm{h\highlight{_e}^\alpha \avg{\norm{w_h}^{\gamma-1}}}_{\leb{q}({\E\cup\partial\W})}
        \Norm{h\highlight{_e}^{-\alpha-1} \jump{{w_h}}}_{\leb{p}({\E\cup\partial\W})},
      \end{split}
    \end{equation}
    which in view of (\ref{eq:bound-on-avg}) gives
    \begin{equation}
      \label{eq:bound-on-jump}
        \int_{\E\cup\partial\W} h\highlight{_e}^{-1} \norm{\jump{\norm{w_h}^{\gamma}}}
        \leq
        C   
        \qp{
          \int_\W \norm{w_h}^{q\qp{\gamma-1}}
        }^{\frac{1}{q}}
        \Norm{h\highlight{_e}^{-\frac{1}{q}-1} \jump{{w_h}}}_{\leb{p}({\E\cup\partial\W})}
    \end{equation}
    again where $\alpha = \frac{1}{q}$. 

    Collecting the three bounds (\ref{eq:bound-on-lap}),
    (\ref{eq:final-bound-on-jump-grad}) and (\ref{eq:bound-on-jump}) and
    substituting into (\ref{eq:step-1}) shows
    \begin{equation}
      \label{eq:bound-step-1}
      \begin{split}
        \qp{\int_\W \norm{w_h}^{\frac{\gamma d}{d-1}}}^{\frac{d-1}{d}} 
        &\leq 
        \qp{\int_\W \norm{w_h}^{q\qp{\gamma - 1}}}^{\frac{1}{q}}
        \bigg(
        \Norm{\Hess_h w_h}_{\leb{p}(\W)}
        +
        \Norm{h\highlight{_e}^{-{\frac{1}{q}}} \jump{\nabla_h{w_h}}}_{\leb{p}({\E\cup\partial\W})}
        \\
        &
        \qquad
        \qquad
        \qquad
        \qquad
        \qquad
        \qquad
        +
        \Norm{h\highlight{_e}^{-\frac{1}{q}-1} \jump{{w_h}}}_{\leb{p}({\E\cup\partial\W})}
        \bigg).
      \end{split}
    \end{equation}
    The main idea of the proof is to now choose $\gamma$ such that
    $\frac{\gamma d}{d-1} = q\qp{\gamma - 1}$. Hence $\gamma =
    \frac{p\qp{d-1}}{d-p}$. Using this and dividing through by the
    first term on the right hand side of (\ref{eq:bound-step-1})
    yields
    \begin{equation}
      \begin{split}
      \qp{\int_\W \norm{w_h}^{\frac{p d}{d-p}}}^{\frac{d-1}{d} - \frac{1}{q}} 
      &\leq
      \bigg(
        \Norm{\Hess_h w_h}_{\leb{p}(\W)}
        +
        \Norm{h\highlight{_e}^{-{\frac{1}{q}}} \jump{\nabla_h{w_h}}}_{\leb{p}({\E\cup\partial\W})}
       \\
       &
       \qquad\qquad \qquad\qquad \qquad\qquad +
        \Norm{h\highlight{_e}^{-\frac{1}{q}-1} \jump{{w_h}}}_{\leb{p}({\E\cup\partial\W})}
        \bigg).
      \end{split}
    \end{equation}
    Now noting that
    \begin{gather}
      \frac{d-1}{d}-\frac{1}{q} = \frac{d-p}{dp}
      \\
      h\highlight{_e}^{-\frac{p}{q}} = h\highlight{_e}^{1-p} \AND
      \\
      h\highlight{_e}^{-\frac{p}{q}-p}= h\highlight{_e}^{1-2p}
    \end{gather}
    yields
    \begin{equation}
      \Norm{w_h}_{\leb{{p^*}}(\W)}
      \leq
      \enorm{w_h}{p}
    \end{equation}
    where $p^* = \frac{pd}{p-d}$ is the \emph{Sobolev conjugate} of
    $p$. This yields the desired result since $p^* > p$ for
    $p\in (1,d)$ and hence we may use the embedding
    ${\leb{p^*}(\W)}\subset\subset {\leb{p}(\W)}$.

    Step 2. For the case $p\in[d, \infty)$ we set $r =
      \frac{dp}{d+p}$. We note that $r < d$ and that the Sobolev
      conjugate of $r$, $r^* = \frac{dr}{d-r} > r$. Following the
      arguments given in Step 1 we arrive at
      \begin{equation}
        \Norm{w_h}_{\leb{{r^*}}(\W)}
        \leq
        \enorm{w_h}{r}.
      \end{equation}
      Note that
      \begin{equation}
        r^* = \frac{rd}{d-r} = \frac{\frac{d^2 p}{d+p}}{d - \frac{dp}{d+p}} = p.
      \end{equation}
      Hence we see that
      \begin{equation}
        \Norm{w_h}_{\leb{p}(\W)}
        =
        \Norm{w_h}_{\leb{r^*}(\W)}
        \leq
        C
        \enorm{w_h}{r}
        \leq 
        C
        \enorm{w_h}{p},
      \end{equation}
      where the final bound follows from Lemma
      \ref{eq:equiv-of-dg-norm}, concluding the proof.
\end{Proof}

\begin{Hyp}[approximability of the finite element space]
  \label{ass:approx}
  Henceforth we will assume the finite element space $\fes$ is chosen
  such that the $\leb{2}(\W)$ orthogonal projection operator satisfies:
  \begin{gather}
    \lim_{h\to 0}\Norm{v - \ltwoproj{}v}_{\leb{p}(\W)} = 0
    \\
    \lim_{h\to 0}\Norm{\nabla v - \nabla_h \qp{\ltwoproj{}v}}_{\leb{p}(\W)} = 0\AND
    \\
    \lim_{h\to 0}\enorm{v - \ltwoproj{}v}{p} = 0.
  \end{gather}
  A choice of $k \geq 2$ satisfies these assumptions.
\end{Hyp}

\begin{The}[stability]  
  \label{the:stability}
  Let $\H[\cdot]$ be defined as in Example \ref{ex:IP-Hess} then the
  dG Hessian is stable in the sense that
  \begin{equation}
    \label{eq:bound-the-liftings}
    \begin{split}
      \Norm{\Hess_h v_h - \H[v_h]}^p_{\leb{p}(\W)^{d\times d}} 
      &\leq
      C\qp{\Norm{l_1[v_h] + l_2[v_h]}^p_{\leb{p}(\W)^{d\times d}}}
      \\
      &\leq
      C\qp{\int_{\E\cup\partial\W} h\highlight{_e}^{1-p} 
      \norm{\jump{\nabla_h v_h}}^p
      +
      h\highlight{_e}^{1-2p}
      \norm{\jump{v_h}}^p
      }.
    \end{split}
  \end{equation}
  Consequently we have
  \begin{equation}
    \label{eq:stab-of-H}
    \Norm{\H[v_h]}^p_{\leb{p}(\W)^{d\times d}}  
    \leq
    C \enorm{v_h}{p}^p
  \end{equation}
\end{The}
\begin{Proof}
  We begin by bounding each of the lifting operators individually. Let
  $q = \tfrac{p}{p-1}$ then by the definition of the $\leb{p}(\W)$
  norm we have that
  \begin{equation}
    \Norm{l_1[v_h]}_{\leb{p}(\W)}
    =
    \sup_{z\in\leb{q}(\W)}
    \int_\W
    \frac{l_1[v_h] z}{\Norm{z}_{\leb{q}(\W)}}.
  \end{equation}
  Let $\ltwoproj{} : \leb{2}(\W) \to \fes$ denote the orthogonal
  projection operator then using the definition of $l_1[\cdot]$
  \eqref{eq:lift-1} we see
  \begin{equation}
    \label{eq:lifting-bound-1}
    \begin{split}
      \Norm{l_1[v_h]}_{\leb{p}(\W)}
      &=
      \sup_{z\in\leb{q}(\W)}
      \int_\W
      \frac{l_1[v_h] \ltwoproj{}z}{\Norm{z}_{\leb{q}(\W)}}
      \\
      &=
      \sup_{z\in\leb{q}(\W)}
      \int_{\E\cup\partial\W}
      \frac{\jump{v_h} \otimes \avg{\nabla_h\qp{\ltwoproj{}z}}}{\Norm{z}_{\leb{q}(\W)}}
      \\
      &\leq
      d^2 
      \sup_{z\in\leb{q}(\W)}
      \frac{\Norm{h\highlight{_e}^{-\alpha} \jump{v_h}}_{\leb{p}({\E\cup\partial\W})} 
        \Norm{\avg{h\highlight{_e}^{\alpha} \nabla_h\qp{\ltwoproj{}z}}}_{\leb{q}({\E\cup\partial\W})}}
             {\Norm{z}_{\leb{q}(\W)}}
      \\
      &\leq 
      d^2 
      \sup_{z\in\leb{q}(\W)}
      \frac{\qp{\Norm{h\highlight{_e}^{-\alpha} \jump{v_h}}_{\leb{p}({\E\cup\partial\W})}^p}^{1/p} 
        \qp{\Norm{\avg{h\highlight{_e}^{\alpha} \nabla_h\qp{\ltwoproj{}z}}}_{\leb{q}({\E\cup\partial\W})}^q}^{1/q}}
             {\Norm{z}_{\leb{q}(\W)}}
    \end{split}
  \end{equation}
  using a H\"older inequality, followed by a discrete H\"older
  inequality and where $\alpha\in\reals$ is some parameter to be
  chosen. 

  Using the definition of the average operator we see
  \begin{equation}
    \Norm{\avg{h\highlight{_e}^{\alpha} \nabla_h\qp{\ltwoproj{}z}}}_{\leb{q}({\E\cup\partial\W})}^q
    \leq
    \tfrac{1}{2} \sum_{K\in\T{}}
    \Norm{h\highlight{_e}^\alpha \nabla_h \qp{\ltwoproj{}z}}_{\leb{q}(\partial K)}^q.
  \end{equation}
  Now using the trace inequality given in Proposition \ref{pro:trace}
  we have
  \begin{equation}
    \Norm{\avg{h\highlight{_e}^{\alpha} \nabla_h\qp{\ltwoproj{}z}}}_{\leb{q}({\E\cup\partial\W})}^q
    \leq
    C 
    \sum_{K\in\T{}}
    h^{q\alpha - 1}
    \Norm{\nabla_h \qp{\ltwoproj{}z}}_{\leb{q}(K)}^q.
  \end{equation}
  Making use of the inverse inequality given in
  Proposition \ref{pro:inverse} we see
  \begin{equation}
    \label{eq:bound-on-l1}
    \Norm{\avg{h\highlight{_e}^{\alpha} \nabla_h\qp{\ltwoproj{}z}}}_{\leb{q}({\E\cup\partial\W})}^q
    \leq
    C 
    \sum_{K\in\T{}}
    h^{q\alpha - 1 - q}
    \Norm{{\ltwoproj{}z}}_{\leb{q}(K)}^q.
  \end{equation}
  We choose $\alpha = 2 - \tfrac{1}{p}$ such that the exponent of $h$
  in the final term of (\ref{eq:bound-on-l1}) is zero. Substituting
  this bound into (\ref{eq:bound-on-l1}) and making use of the
  stability of the $\leb{2}(\W)$ orthogonal projection in
  $\leb{p}(\W)$ \highlight{\cite{CrouzeixThomee:1987} \margnote{Ref A, comment (13)}} we see that
  \begin{equation}
    \label{eq:result-on-l1}
    \begin{split}
      \Norm{l_1[v_h]}_{\leb{p}(\W)}^p
      &\leq
      C \Norm{h\highlight{_e}^{\tfrac{1}{p}-2} \jump{v_h}}_{\leb{p}({\E\cup\partial\W})}^p
      \\
      &\leq
      C h\highlight{_e}^{1-2p}\Norm{\jump{v_h}}_{\leb{p}({\E\cup\partial\W})}^p.
    \end{split}
  \end{equation}

  The bound on $l_2[\cdot]$ is achieved using much the same
  argument. Following the steps given in (\ref{eq:lifting-bound-1}) it
  can be verified that
  \begin{equation}
    \label{eq:lifting-bound-2}
    \Norm{l_2[v_h]}_{\leb{p}(\W)}
    \leq 
    d^2 
    \sup_{z\in\leb{q}(\W)}
    \frac{\qp{\Norm{h^{-\beta} \jump{\nabla_h v_h}}_{\leb{p}({\E\cup\partial\W})}^p}^{1/p} 
      \qp{\Norm{\avg{h^{\beta} {\ltwoproj{}z}}}_{\leb{q}({\E\cup\partial\W})}^q}^{1/q}}
         {\Norm{z}_{\leb{q}(\W)}}
  \end{equation}
  for some $\beta\in\reals$. To bound the average term, we follow the
  same steps (without the inverse inequality)
  \begin{equation}
    \begin{split}
      \Norm{\avg{h\highlight{_e}^{\beta} {\ltwoproj{}z}}}_{\leb{q}({\E\cup\partial\W})}^q
      &\leq
      \tfrac{1}{2}
      \sum_{K\in\T{}}
      \Norm{{h^{\beta} {\ltwoproj{}z}}}_{\leb{q}(\partial K)}^q
      \\
      &\leq
      C
      \sum_{K\in\T{}}
      h^{q\beta -1}\Norm{{{\ltwoproj{}z}}}_{\leb{q}(K)}^q.
    \end{split}
  \end{equation}
  We choose $\beta = 1-\tfrac{1}{p}$ such that the exponent of $h$
  vanishes and substitute into (\ref{eq:lifting-bound-2}) to find
  \begin{equation}
    \label{eq:result-on-l2}
    \begin{split}
      \Norm{l_2[v_h]}_{\leb{p}(\W)}^p
      &\leq
      C \Norm{h\highlight{_e}^{\tfrac{1}{p}-1} \jump{v_h}}_{\leb{p}({\E\cup\partial\W})}^p
      \\
      &\leq
      C h\highlight{_e}^{1-p}\Norm{\jump{v_h}}_{\leb{p}({\E\cup\partial\W})}^p.
    \end{split}
  \end{equation}
  The result (\ref{eq:bound-the-liftings}) follows noting the
  definition of $\H$ given in (\ref{eq:definition-of-hessian}), a
  Minkowski inequality and the two results (\ref{eq:result-on-l1}) and
  (\ref{eq:result-on-l2}).

  To see (\ref{eq:stab-of-H}) it suffices to again use a Minkowski
  inequality, together with (\ref{eq:definition-of-hessian}) and the
  two results (\ref{eq:result-on-l1}) and (\ref{eq:result-on-l2}).
\end{Proof}

\highlight{
  \begin{Cor}[strong convergence of the dG-Hessian]
    \label{cor-strong-conv-of-dg-hess}
  Given a smooth $v\in\cont{\infty}_0(\W)$, with
  $\ltwoproj:\leb{2}(\W)\to\fes$ being the $\leb{2}$ orthogonal
  projection operator we have that
  \begin{equation}
    \Norm{\Hess v - \H[\ltwoproj{v}]}_{\leb{p}(\W)^{d\times d}}
    \leq
    C \enorm{v - \ltwoproj{}v}{p}.
  \end{equation}
  Hence using the approximation properties given in Assumption
  \ref{ass:approx}, we have that $\H[\ltwoproj v] \to \Hess v$
  strongly in $\leb{p}(\W)^{d\times d}$.
  \margnote{Ref A, comment (4)}
\end{Cor}
}
\subsection{Numerical minimisation problem and discrete Euler--Lagrange equations}

The properties of the IP-Hessian allow us to define the following
numerical scheme: To seek $u_h\in\fes$ such that
\begin{equation}
  \label{eq:discrete-min}
  \cJ_h[u_h; p] = \inf_{v_h\in\fes} \cJ_h[v_h;p].
\end{equation}
Let $\cD[v_h]:=\trace{\H[v_h]}$ then the discrete action functional
$\cJ_h$ is given by
\begin{equation}
  \cJ_h[v_h; p] := \int_\W \frac{1}{p}\norm{\cD[v_h]}^p + f v_h
  +
  \frac{\sigma}{p} \qp{
  \int_{\E\cup \partial\W} h\highlight{_e}^{1-p}\norm{\jump{\nabla_h v_h}}^p
  +
  h\highlight{_e}^{1-2p}\norm{\jump{v_h}}^p
  }
\end{equation}
where $\sigma>0$ is a \emph{penalisation parameter}.

Let
\begin{equation}
  \label{eq:discrete-bilinear}
  \begin{split}
    \bih{u_h}{\Phi}{p}
    &:=
    \int_\W
    \norm{\cD[u_h]}^{p-2}\cD[u_h] \cD[\Phi]
    \\
    &\qquad
    +
    \sigma
    \bigg(
      \int_{\E\cup \partial\W} 
      h\highlight{_e}^{1-p}
      \norm{\jump{\nabla_h u_h}}^{p-2}
      \jump{\nabla_h u_h}
      \jump{\nabla_h \Phi}
      \\
      &\qquad\qquad\qquad\qquad\qquad\qquad+
      h\highlight{_e}^{1-2p}\norm{\jump{v_h}}^{p-2}
      \jump{u_h}
      \jump{\Phi}
      \bigg)
  \end{split}
\end{equation}

The associated (weak) discrete Euler--Lagrange equations to the
problem are to seek $\qp{u_h, \H[u_h]} \in \fes\times\fes^{d\times d}$
such that
\begin{equation}
  \label{eq:discrete-EL-eqns}
  \bih{u_h}{\Phi}{p}
  = 
  \int_\W f \Phi
  \Foreach \Phi \in\fes,
\end{equation}
where $\H$ is defined in Example \ref{ex:IP-Hess}.

\begin{The}[coercivity]
  \label{the:coer-and-cont}
  Let $f\in\leb{q}(\W)$ and $\{u_h, \fes\}$ be the finite element
  sequence satisfying the discrete minimisation problem
  (\ref{eq:discrete-min}) then we have that there exists constants
  $C\highlight{=C(p)}>0$\margnote{Ref A, comment (5)} and $\gamma\geq 0$ such that
  \begin{equation}
    \label{eq:coercive-discrete-functional}
    \cJ_h[u_h;p] \geq C \enorm{u_h}{p}^p - \gamma.
  \end{equation}
  Equivalently let $\bih{\cdot}{\cdot}{p}$ be defined as in
  (\ref{eq:discrete-bilinear}) then 
  \begin{equation}
    \label{eq:coercive-discrete-bilinear}
    \bih{u_h}{u_h}{p} \geq C\enorm{u_h}{p}^p.
  \end{equation}
\end{The}
\begin{Proof}
  We have by definition of $\enorm{\cdot}{p}$ that
  \begin{equation}
    \enorm{u_h}{p}^p
    = 
    \Norm{\Hess_h u_h}^p_{\leb{p}(\W)}
    +
    h\highlight{_e}^{1-p}\Norm{\jump{\nabla_h u_h}}^p_{\leb{p}({\E\cup\partial\W})}
    +
    h\highlight{_e}^{1-2p}\Norm{\jump{u_h}}^p_{\leb{p}({\E\cup\partial\W})}.
  \end{equation}
  We see by a Minkowski inequality that
  \begin{equation}
    \begin{split}
      \enorm{u_h}{p}^p 
      &\leq
      \Norm{\Hess_h u_h - \H[u_h]}^p_{\leb{p}(\W)}
      +
      \Norm{\H[u_h]}^p_{\leb{p}(\W)}
      \\
      &
      \qquad
      +
      h\highlight{_e}^{1-p}\Norm{\jump{\nabla_h u_h}}^p_{\leb{p}({\E\cup\partial\W})}
      +
      h\highlight{_e}^{1-2p}\Norm{\jump{u_h}}^p_{\leb{p}({\E\cup\partial\W})}. 
    \end{split}
  \end{equation}
  Hence, using the stability of the discrete Hessian given in Theorem
  \ref{the:stability} we have that
  \begin{equation}
    \begin{split}
      \enorm{u_h}{p}^p 
      &\leq
      \Norm{\H[u_h]}^p_{\leb{p}(\W)}
      +
      \qp{1+C\highlight{(p)}}\bigg(
        h\highlight{_e}^{1-p}\Norm{\jump{\nabla_h u_h}}^p_{\leb{p}({\E\cup\partial\W})}
        \\
        &\qquad\qquad\qquad\qquad\qquad\qquad +
        h\highlight{_e}^{1-2p}\Norm{\jump{u_h}}^p_{\leb{p}({\E\cup\partial\W})} 
      \bigg)
      \\
      &\leq
      C(p) \bih{u_h}{u_h}{p},
    \end{split}
  \end{equation}
  where we have made use of a piecewise equivalent of Proposition
  \ref{prop:equiv-of-norms} hence showing
  (\ref{eq:coercive-discrete-bilinear}). The result
  (\ref{eq:coercive-discrete-functional}) follows using a similar
  argument.
\end{Proof}

\begin{Lem}[relative compactness]
  \label{lem:rel-compact}
  Let $\{ v_h, \fes \}$ be a finite element sequence that is bounded
  in the $\enorm{\cdot}{p}$ norm. Then the sequence is relatively
  compact in $\leb{p}(\W)$.
\end{Lem}
\begin{Proof}
  The proof is an application of Kolmogorov's Compactness Theorem
  noting the result of Lemma \ref{lem:discrete-sob-embed} which infers
  boundedness of the finite element sequence in $\leb{p}(\W)$.
\end{Proof}

\begin{Lem}[limit]
  \label{lem:limit}
  Given a finite element sequence $\{ v_h, \fes \}$ that is bounded
  in the $\enorm{\cdot}{p}$ norm, there exists a function
  $v\in\sobz{2}{p}(\W)$ such that as $h\to 0$ we have, up to a
  subsequence, $v_h \rightharpoonup v$ weakly in $\leb{p}(\W)$. Moreover, $\H[v_h]
  \rightharpoonup \Hess v$ weakly in $\leb{p}(\W)^{d\times d}$.
\end{Lem}
\begin{Proof}
  Lemma \ref{lem:rel-compact} infers that we may find a
  $v\in\leb{p}(\W)$ which is the limit of our finite element
  sequence. To prove that $v\in\sobz{2}{p}(\W)$ we must show that our
  sequence of discrete Hessians converge to $\Hess v$.

  Recall Theorem \ref{the:stability} gave us that
  \begin{equation}
    \Norm{\H[v_h]}_{\leb{p}(\W)^{d\times d}}  
    \leq
    C \enorm{v_h}{p}.
  \end{equation}
  As such, we may infer the (matrix valued) finite element sequence
  $\{ \H[v_h], \fes^{d\times d} \}$ is bounded in ${\leb{p}(\W)^{d\times d}}$. Hence
  we have that $ \H[v_h] \rightharpoonup \geomat X \in \leb{p}(\W)^{d\times
    d}$ weakly for some matrix valued function $\geomat X$.

  Now we must verify that $\geomat X = \Hess v$. For each
  $\phi\in\cont{\infty}_0(\W)$ we have that
  \begin{equation}
    \int_\W
    \H[v_h] {\ltwoproj \phi}
    =
    \int_\W \Hess_h v_h {\ltwoproj \phi}
    -
    \int_{\E} \tjump{\nabla_h v_h} \avg{{\ltwoproj \phi}}
    +
    \int_{\E\cup\partial\W} \jump{v_h} \otimes \avg{\nabla_h \qp{\ltwoproj \phi}}.
  \end{equation}
  Note that
  \highlight{
  \begin{equation}
    \begin{split}
      \int_\W \Hess_h v_h {\ltwoproj \phi}
      &= 
      -\int_\W \nabla_h v_h \otimes \nabla_h \qp{{\ltwoproj \phi}}
      +
      \int_{\E} \tjump{\nabla_h v_h} \avg{{\ltwoproj \phi}}
      \\
      &\qquad +
      \int_{\E\cup\partial\W} \jump{\ltwoproj \phi} \otimes \avg{\nabla_h v_h}
      \\
      &=
      \int_\W v_h \Hess_h \qp{{\ltwoproj \phi}}
      +
      \int_{\E} \tjump{\nabla_h v_h} \avg{{\ltwoproj \phi}}
      - 
      \tjump{\nabla_h \qp{\ltwoproj\phi}} \avg{{v_h}}
      \\
      &
      \qquad 
      +
      \int_{\E\cup\partial\W} \jump{\ltwoproj \phi} \otimes \avg{\nabla_h v_h}
      -
      \jump{v_h} \otimes \avg{\nabla_h \qp{\ltwoproj \phi}}
      \\
      &=
      \int_\W v_h \H[\ltwoproj \phi]
      +
      \int_{\E} \tjump{\nabla_h v_h} \avg{{\ltwoproj \phi}}
      -
      \int_{\E\cup\partial\W}
      \jump{v_h} \otimes \avg{\nabla_h \qp{\ltwoproj \phi}}
    \end{split}
  \end{equation}
  \margnote{Ref A, comment (4)}}
  As such, we have that
  \begin{equation}
    \begin{split}
      \int_\W 
      \geomat X \phi
      &=
      \lim_{h\to 0}
      \int_\W 
      \H[v_h]{\ltwoproj \phi}
      \\
      &=
      \lim_{h\to 0}
      \highlight{\int_\W v_h \H[\ltwoproj \phi]}
      \\
      &= 
      \int_\W v \Hess \phi
    \end{split}
  \end{equation}
  \highlight{by the strong convergence of the dG Hessian in Corollary
    \ref{cor-strong-conv-of-dg-hess}}. Hence we
  have that $\geomat X = \Hess v$ in the distributional sense.
\end{Proof}

\begin{Lem}[apriori bound]
  \label{lem:apriori}
  Let $f\in\leb{q}(\W)$, with $q=\tfrac{p}{p-1}$ and let $\{ u_h, \fes
  \}$ be the finite element sequence satisfying
  (\ref{eq:discrete-min}), then we have the following apriori bound:
  \begin{equation}
    \enorm{u_h}{p} \leq \qp{C \Norm{f}_{\leb{q}(\W)}}^{q/p}.
  \end{equation}
\end{Lem}
\begin{Proof}
  Using the coercivity condition given in Theorem
  \ref{the:coer-and-cont} and the definition of the weak
  Euler--Lagrange equations we have
  \begin{equation}
    \begin{split}
      \enorm{u_h}{p}^p 
      &\leq 
      C \bih{u_h}{u_h}{p}
      \\
      &\leq
      C \int_\W f u_h.
    \end{split}
  \end{equation}
  Now using a H\"older inequality and the discrete Sobolev embedding
  given in Lemma \ref{lem:discrete-sob-embed} we see
  \begin{equation}
    \begin{split}
      \enorm{u_h}{p}^p 
      &\leq 
      C
      \Norm{f}_{\leb{q}(\W)}
      \Norm{u_h}_{\leb{p}(\W)}
      \\
      &\leq 
      C
      \Norm{f}_{\leb{q}(\W)}
      \enorm{u_h}{p}.
    \end{split}
  \end{equation}
  Upon simplifying, we obtain the desired result.
\end{Proof}

\begin{The}[convergence]
  Let $f\in\leb{q}(\W)$, with $q=\tfrac{p}{p-1}$ and suppose $\{ u_h,
  \fes \}$ is the finite element sequence generated by solving the
  nonlinear system (\ref{eq:discrete-EL-eqns}), then we have that
  \begin{itemize}
  \item $u_h\to u$ in $\leb{p}(\W)$ and
  \item $\H[u_h] \to \Hess u$ in $\leb{p}(\W)^{d\times d}$.
  \end{itemize}  
  where $u\in\sobz{2}{p}(\W)$ be the unique solution to the
  $p$--biharmonic problem \eqref{eq:p-biharm}.
\end{The}
\begin{Proof}
  Given $f\in\leb{q}(\W)$ we have that, in view of Lemma
  \ref{lem:apriori}, the finite element sequence $\{ u_h, \fes\}$ is
  bounded in the $\enorm{\cdot}{p}$ norm. As such we may apply Lemma
  \ref{lem:limit} which shows that there exists a (weak) limit to the
  finite element sequence $\{u_h,\fes \}$ which we shall call $u^\ast$. We
  must now show that $u^\ast = u$, the solution of the $p$--biharmonic
  problem.
  
  By Corollary \ref{cor:semicont} $\cJ[\cdot]$ is weakly lower
  semicontinuous, hence we have that
  \begin{equation}
    \begin{split}
      \cJ[u^\ast] 
      &\leq 
      \liminf_{h\to 0} 
      \qb{
        \frac{1}{p}\Norm{\cD[u_h]}_{\leb{p}(\W)}^p + \int_\W f u_h
      }
      \\
      &\leq
      \liminf_{h\to 0} 
      \bigg[
        \frac{1}{p}\Norm{\cD[u_h]}_{\leb{p}(\W)}^p + \int_\W f u_h
        \\
        &
        \qquad 
        \qquad 
        +
        \frac{\sigma}{p}
        \qp{
        h\highlight{_e}^{1-p}
        \Norm{\jump{\nabla_h u_h}}_{\leb{p}(\W)}^p 
        +
        h\highlight{_e}^{1-2p}
        \Norm{\jump{u_h}}_{\leb{p}(\W)}^p 
        }
      \bigg].
      \\
      &=
      \liminf_{h\to 0} 
      \cJ_h[u_h].
    \end{split}
  \end{equation}
  Now owing to Assumption \ref{ass:approx} we have that for any
  $v\in\cont{\infty}_0(\W)$ that
  \begin{equation}
    \begin{split}
      \cJ[v]
      &=
      \liminf_{h\to 0}
      \bigg[   
        \frac{1}{p}\Norm{\cD[\ltwoproj{}v]}_{\leb{p}(\W)}^p + \int_\W f \ltwoproj{}v
        \\
        &
        \qquad 
        \qquad 
        +
        \frac{\sigma}{p}
        \qp{
          h\highlight{_e}^{1-p}
          \Norm{\jump{\nabla_h \qp{\ltwoproj{}v}}}_{\leb{p}(\W)}^p 
          +
          h\highlight{_e}^{1-2p}
          \Norm{\jump{\ltwoproj{}v}}_{\leb{p}(\W)}^p 
        }
        \bigg]
      \\
      &= 
      \liminf_{h\to 0}
      \cJ_h[\ltwoproj{} v]
    \end{split}
  \end{equation}
  By the definition of the discrete scheme we have that
  \begin{equation}
    \cJ[u^\ast] \leq \cJ_h[u_h] \leq \cJ_h[\ltwoproj{} v] = \cJ[v].
  \end{equation}
  Now, since $v$ was a generic element we may use the density of
  $\cont{\infty}_0(\W)$ in $\sobz{2}{p}(\W)$ and that since $u$ is the
  unique minimiser we must have that $u^\ast = u$.
\end{Proof}

\begin{Rem}[provable rates for the $2$--biharmonic problem]
  \label{rem:provable-rates}
  In the papers \cite{SuliMozolevski:2007, GeorgoulisHouston:2009}
  rates of convergence are given for the $2$--biharmonic problem,
  these are
  \begin{gather}
    \Norm{u-u_h} = \Oh(h^2) \text{ for } k = 2
    \qquad
    \Norm{u-u_h} = \Oh(h^{k+1}) \text{ for } k > 2
    \\
    \enorm{u-u_h}{p} = \Oh(h^{k-1}).
  \end{gather}
  Note that for piecewise quadratic finite elements the convergence
  rate is suboptimal in $\leb{2}(\W)$. 
\end{Rem}

\section{Numerical experiments}
\label{sec:numerics}

In this section we summarise some numerical experiments conducted in
the method presented in Section \ref{sec:discretisation}.

\begin{Rem}[implementation issues]
  The numerical experiments were conducted using the \dolfin interface
  for \fenics \cite{LoggWells:2010}. The graphics were generated using
  \gnuplot and \paraview.
  
  For computational efficiency, we chose to represent $\cD[u_h]$ as an
  auxiliary variable in the mixed formulation, which only requires
  one additional variable, as apposed to the full discrete Hessian
  $\H[u_h]$ which would require $d^2$ (or $\tfrac{d^2 + d}{2}$ if one
  uses symmetry of $\H$). We note that this is only possible due to
  the structure of the problem, \ie that $L = L(\geovec x, u,\nabla
  u,\Delta u)$ and would not be possible in a general setting.
\end{Rem}

\subsection{Benchmarking}
\label{sec:benchmark}

\begin{figure}[ht]
  \caption[]
          {
            \label{fig:benchmark-p-2} 
            Section \ref{sec:benchmark} -- Numerical experiment benchmarking
            the numerical method for the $2$--biharmonic problem. We
            fix $f$ such that the solution $u$ is given by
            (\ref{eq:benchmark-regular}). We plot the log of the error
            together with its estimated order of convergence. We study
            the $\leb{p}(\W)$ norms of the error of the finite element
            solution $u_h$ as well as the represented auxiliary
            variable $\cD[u_h]$ for the dG method
            \eqref{eq:discrete-EL-eqns} with $k = 2,3,4$. We also give
            a solution plot. We observe that the method achieves the rates given in Remark
            \ref{rem:provable-rates}}
            \subfigure[][finite element approximation to (\ref{eq:benchmark-regular}).]
                      {
                        \includegraphics[scale=\figscale, width=0.5\figwidth, angle=180]
                                        {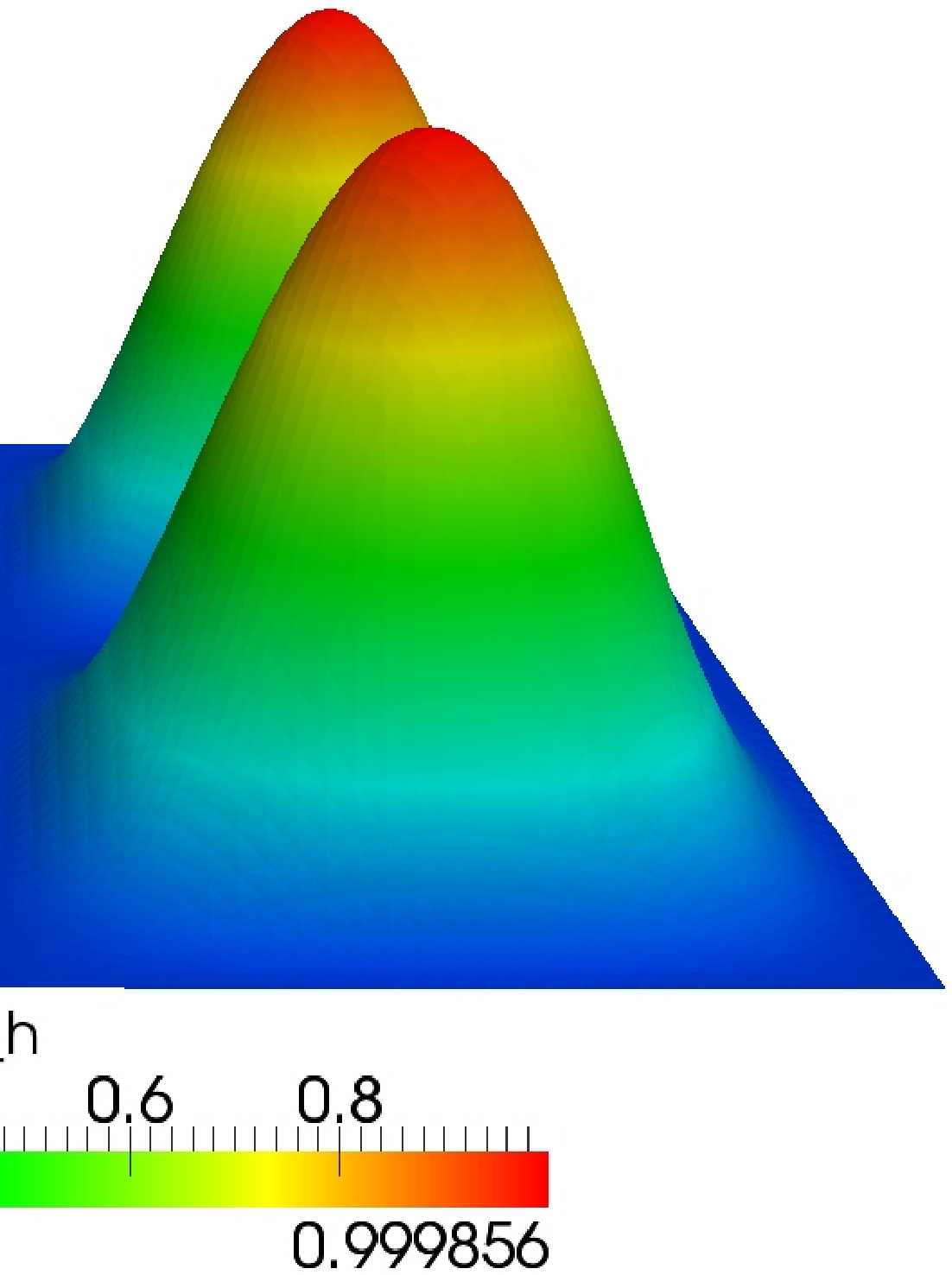} 
                      }
                      \hfill
            \subfigure[][$k = 2$, piecewise quadratic FEs.]
                      {
                        \includegraphics[scale=\figscale, width=0.33\figwidth, angle=270]
                                        {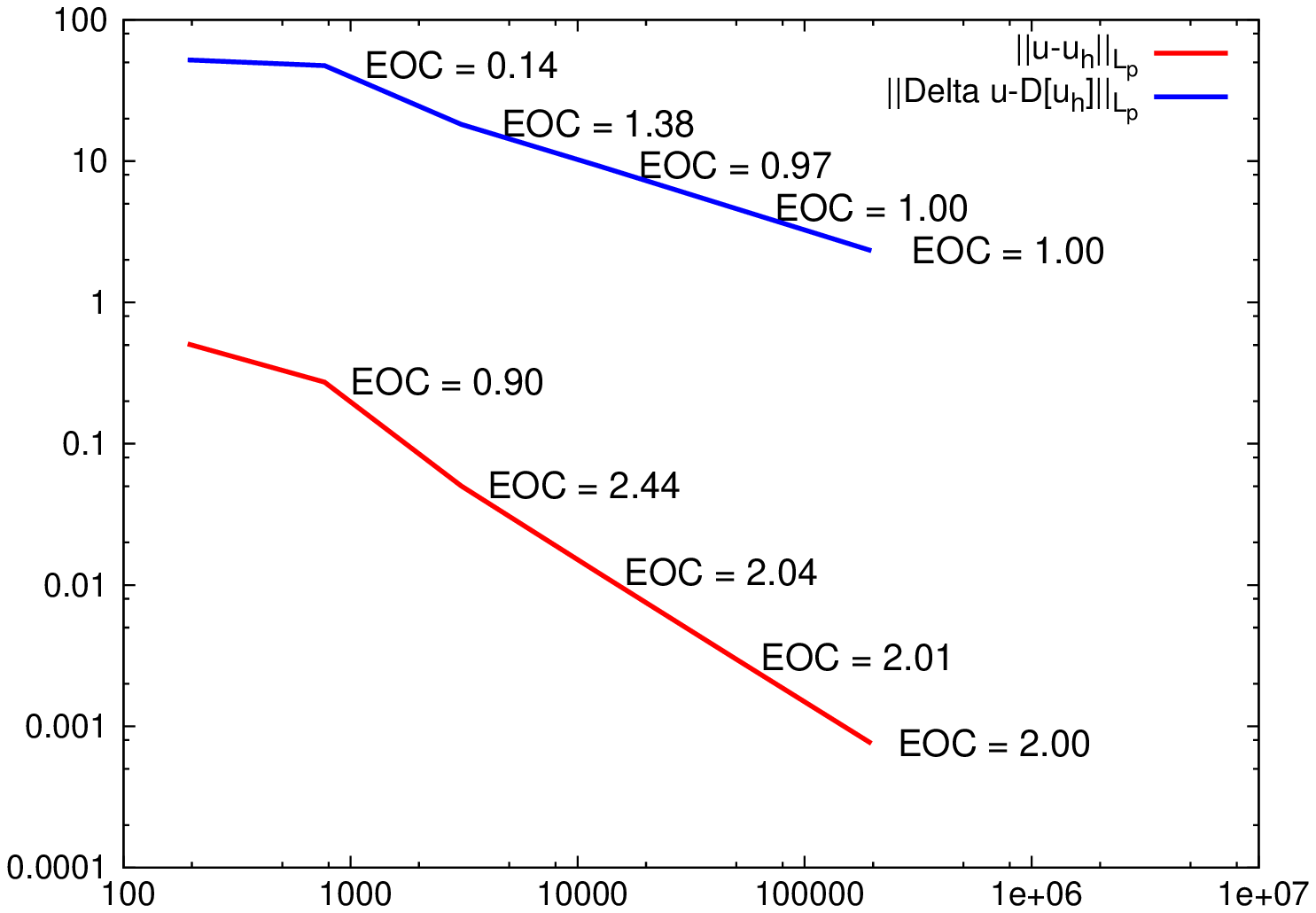}
                      } 
                      \hfill
            \subfigure[][$k = 3$, piecewise cubic FEs.]
                      {
                        \includegraphics[scale=\figscale, width=0.33\figwidth, angle=270]
                                        {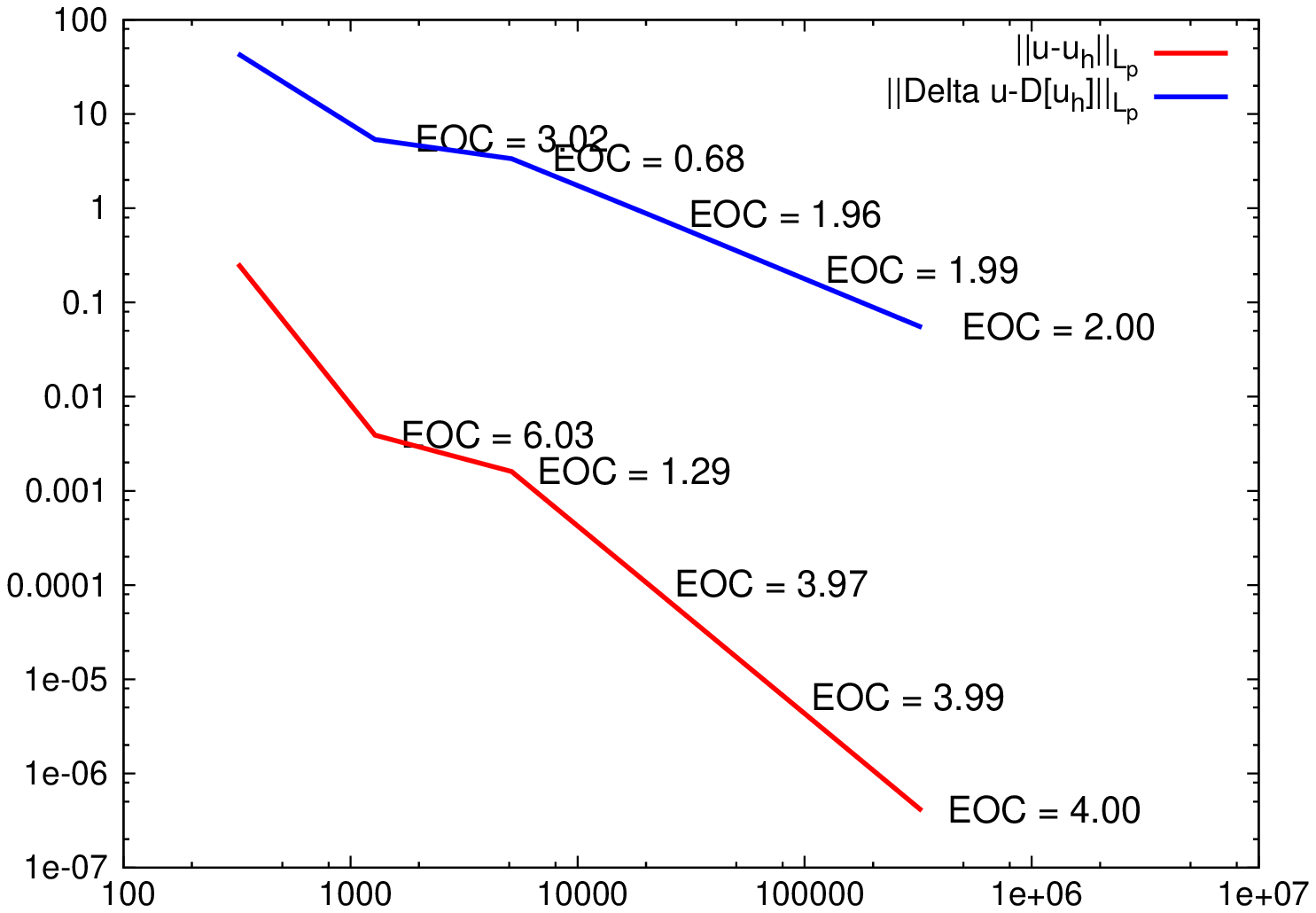}
                      }
                      \hfill
            \subfigure[][$k = 4$, piecewise quartic FEs.]
                      {
                        \includegraphics[scale=\figscale, width=0.33\figwidth, angle=270]
                                        {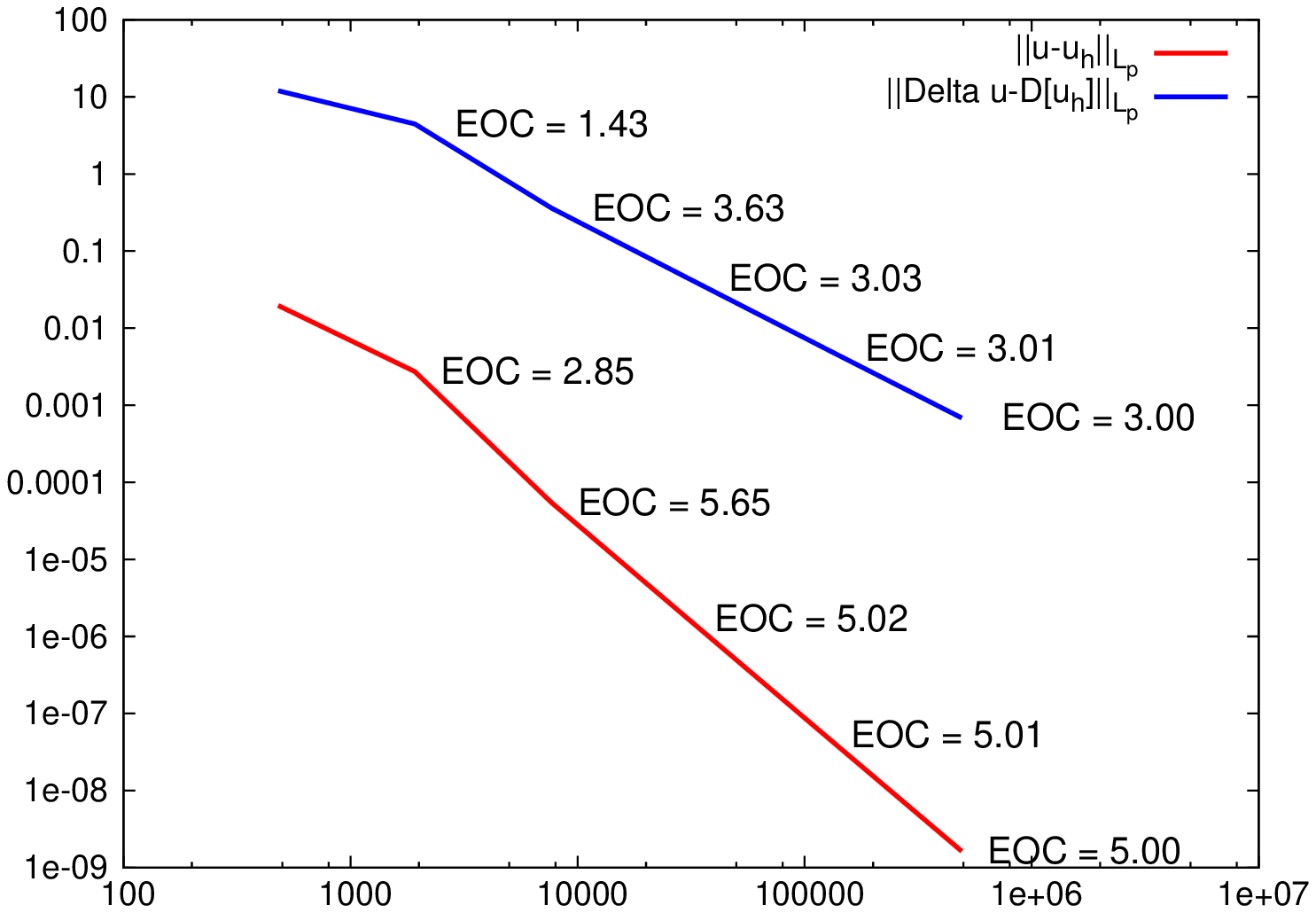}
                      }
\end{figure}

\begin{figure}[ht]
  \caption[]
          {
            \label{fig:benchmark-p-3} 
            Section \ref{sec:benchmark} -- The same test as in Figure
            \ref{fig:benchmark-p-2} for the $2.1$--biharmonic problem,
            \ie $p=2.1$ for $k=2 \AND 3$. }
          \begin{center}
            \subfigure[][$k = 2$, piecewise quadratic FEs.]
                      {
                        \label{fig:benchmark-p-3-2} 
                        \includegraphics[scale=\figscale, width=0.33\figwidth, angle=270]
                                        {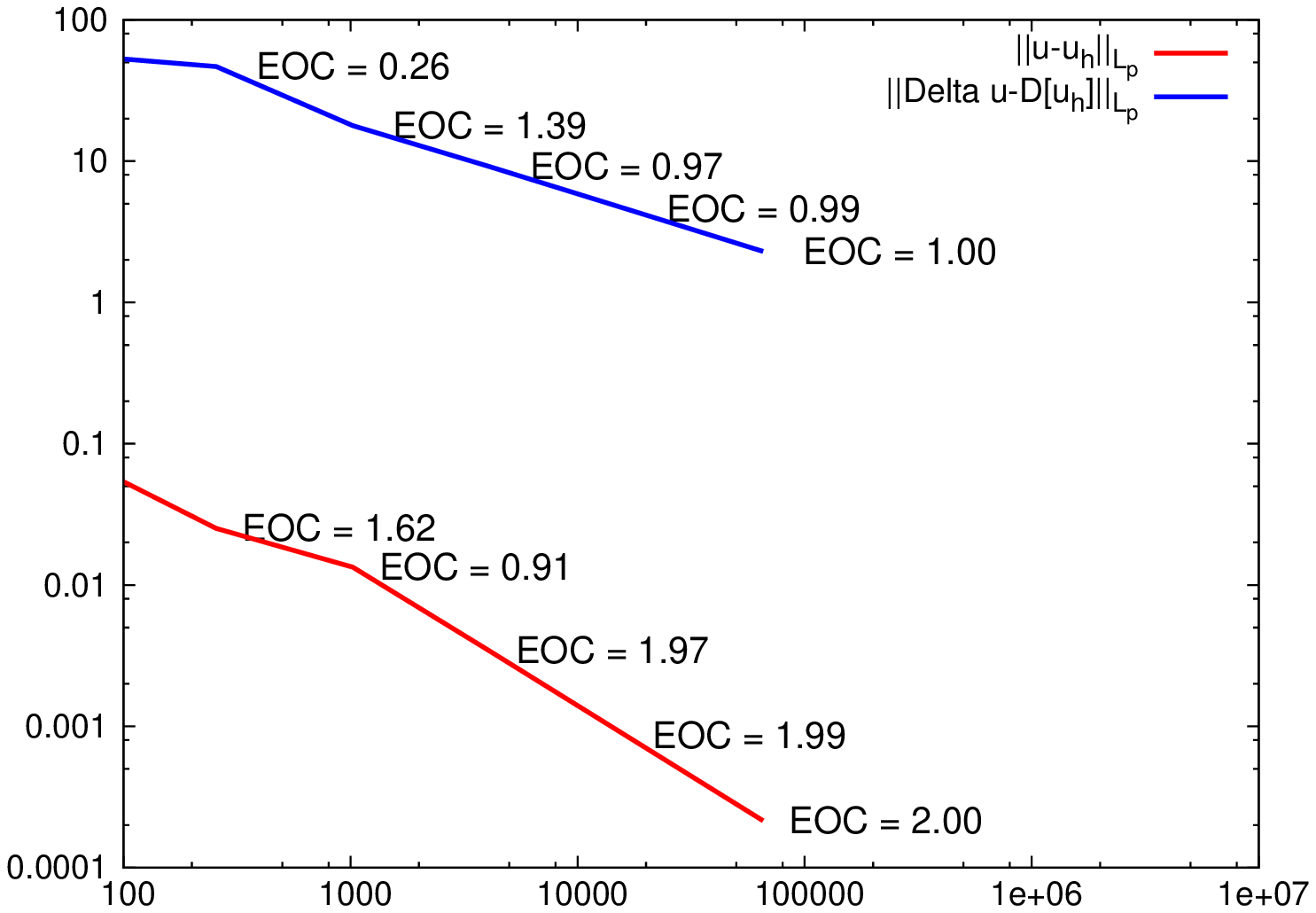}
                      } 
                      \hfill
            \subfigure[][$k = 3$, piecewise cubic FEs.]
                      {
                        \label{fig:benchmark-p-3-3} 
                        \includegraphics[scale=\figscale, width=0.33\figwidth, angle=270]
                                        {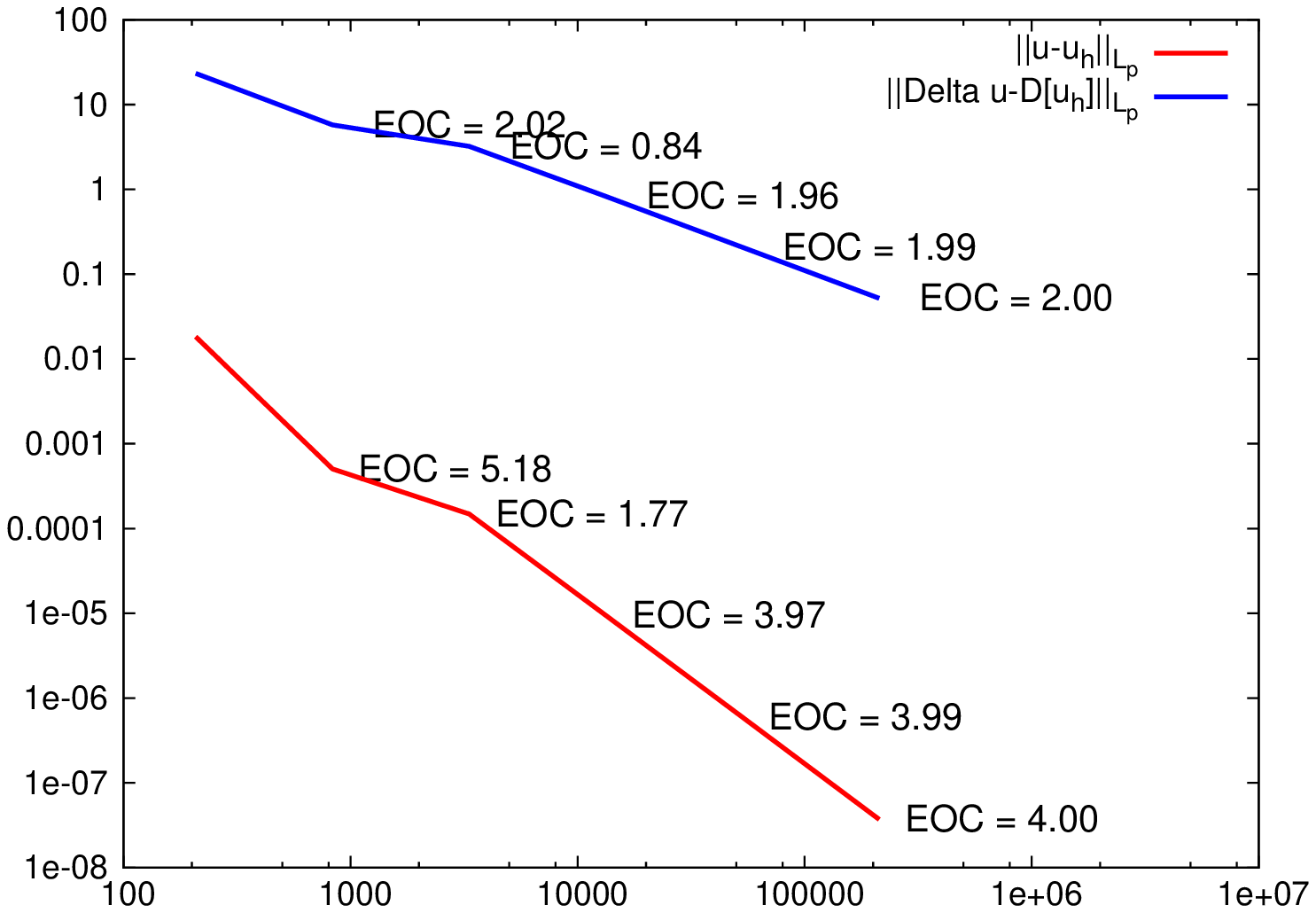}
                      }  
          \end{center}
\end{figure}

\begin{figure}[ht]
  \caption[]
          {
            \label{fig:benchmark-p-4} 
            Section \ref{sec:benchmark} -- The same test as in Figure
            \ref{fig:benchmark-p-3} for the $10$--biharmonic
            problem,\ie $p=10$.}
          \begin{center}
            \subfigure[][$k = 2$, piecewise quadratic FEs.]
                      {
                        \includegraphics[scale=\figscale, width=0.33\figwidth, angle=270]
                                        {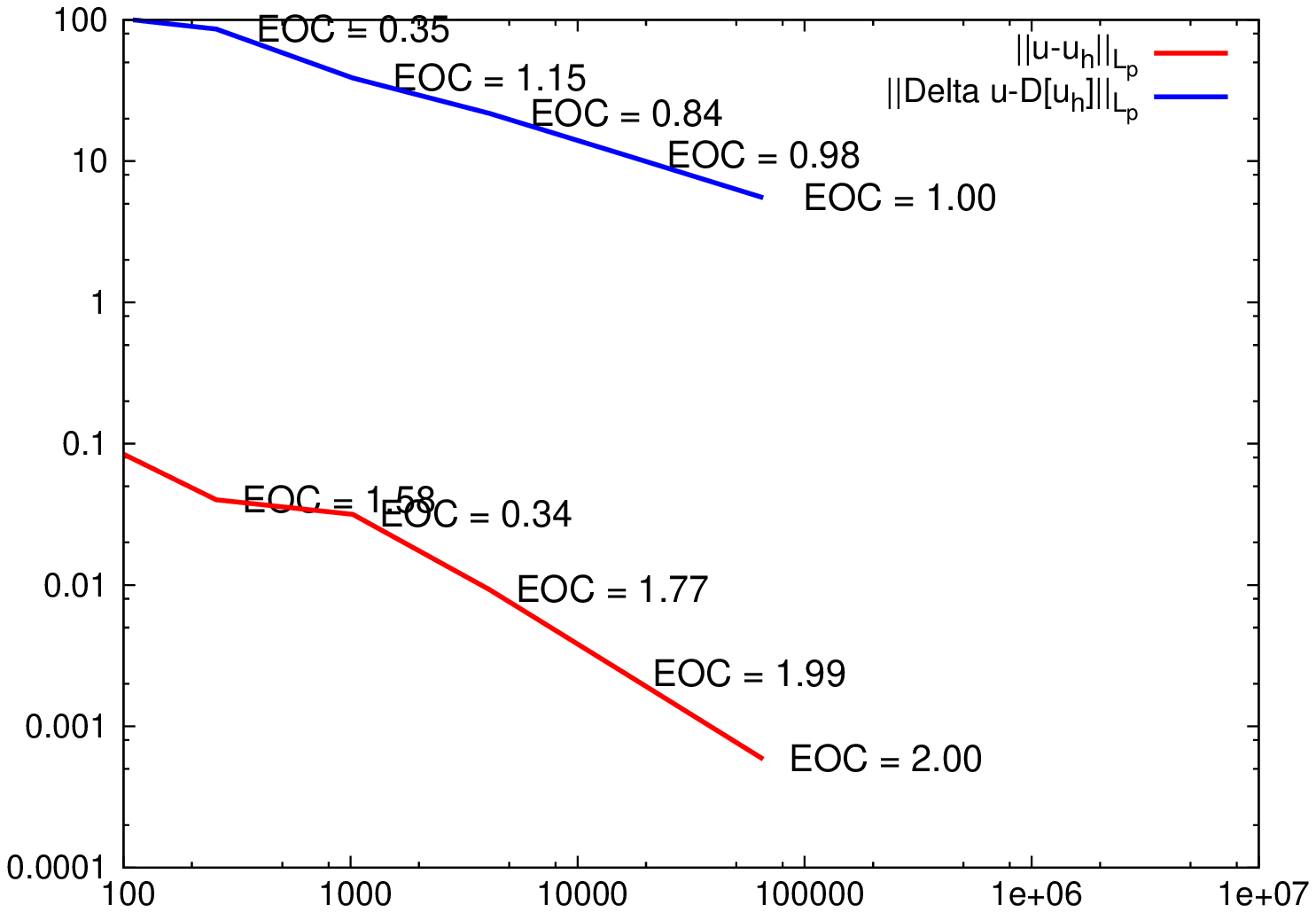}
                      } 
                      \hfill
            \subfigure[][$k = 3$, piecewise cubic FEs.]
                      {
                        \includegraphics[scale=\figscale, width=0.33\figwidth, angle=270]
                                        {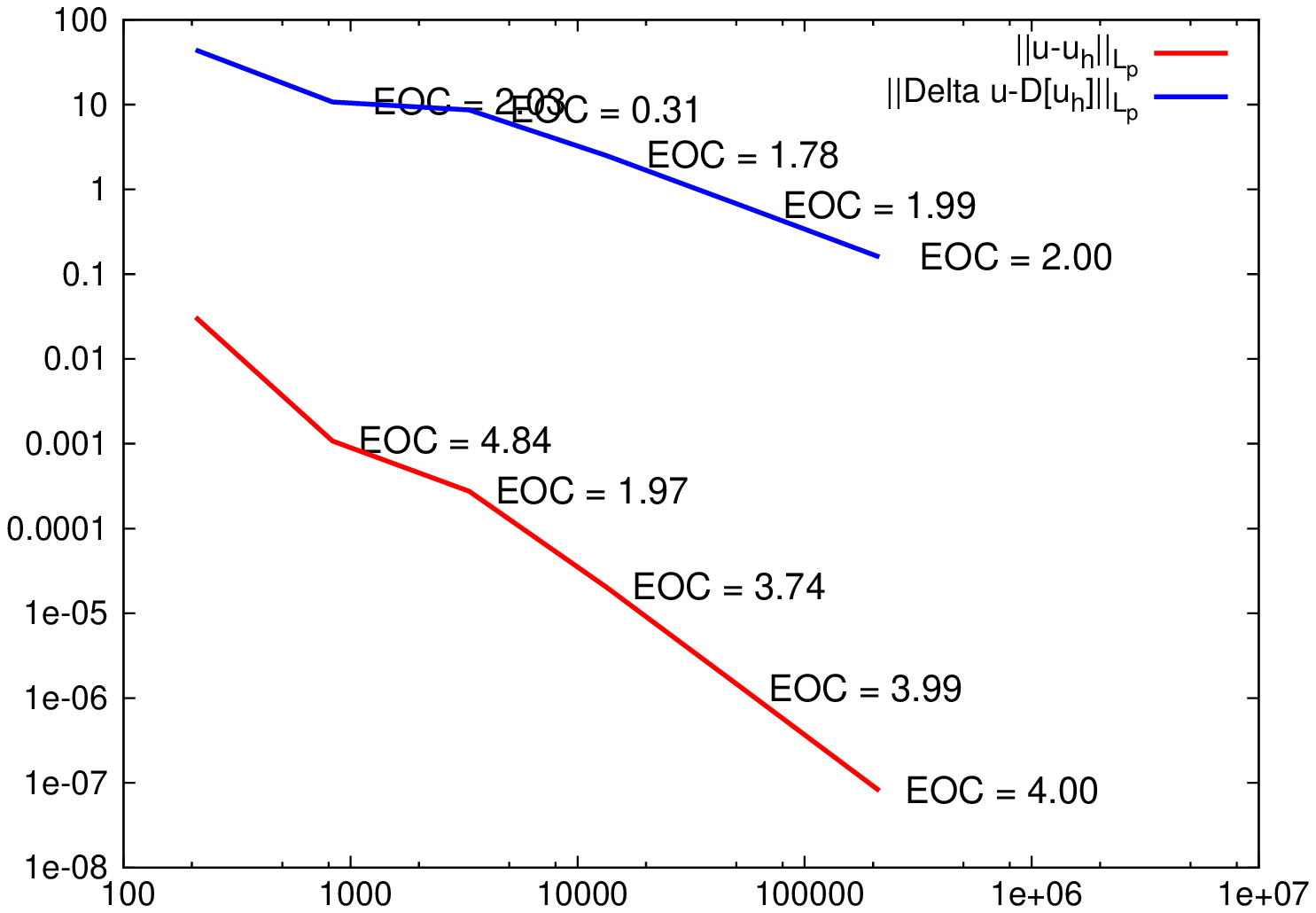}
                      }
                    \end{center}
\end{figure}

The aims of this section are to test the robustness of the numerical
method for a model test solution of the $p$--biharmonic problem. We
show the method achieves the provable rates for $p = 2$ (Figure
\ref{fig:benchmark-p-2}) and numerically gauge the convergence rates
for $p > 2$ (Figures \ref{fig:benchmark-p-3} and
\ref{fig:benchmark-p-4}). To that end,\highlight{we take $\T{}$ to be an unstructured Delaunay triangulation of the square $\W = [0,1]^2$\margnote{Ref A, comment (21)}}. We fix $d=2$, let $\geovec x =
\Transpose{\qp{x,y}}$ and choose $f$ such that
\begin{equation}
  \label{eq:benchmark-regular}
  u(\geovec x) := \sin{2 \pi x}^2\sin{2 \pi y}^2.
\end{equation}
Note that this is comparable to the numerical experiment
\cite[Section 6.1]{GeorgoulisHouston:2009}.

\begin{Rem}[computational observations]
  \label{rem:super}
  Computationally, the convergence rates we observe are that
  \begin{gather}
    \Norm{u-u_h}_{\leb{p}(\W)} = 
    \begin{cases}
      \Oh(h^2) \text{ when $k = 2$}
      \\
      \Oh(h^{k+1}) \text{ otherwise }
    \end{cases}
  \end{gather}
  and that
  \begin{gather}
    \Norm{\Delta u - \cD[u_h]}_{\leb{p}(\W)} =
      \Oh(h^{k-1}).
  \end{gather}
\end{Rem}

\begin{Rem}[representation of $\H$]
  Note that the dG Hessian $\H$ may be represented in a finite element
  space with different degree to $u_h \in \fes$. Let $\fesW :=
  \poly{k-1}(\T{})$, the proof of Theorem
  \ref{the:fully-generalised-fe-hessian} infers that we may choose to
  represent $\H[u_h]\in\fesW^{d\times d}$. For clarity of exposition
  we choose to use $\H[u_h]\in\fes^{d\times d}$, however, we see no
  difficulty extending the arguments presented to the lower degree dG
  Hessian.

  Numerically we observe the same convergence rates as in Remark
  \ref{rem:super} for the lower degree dG Hessian.
\end{Rem}

\section{Conclusion and outlook}

In this work we presented a dG finite element method for the
$p$--biharmonic problem. To do this we introduced a auxiliary
variable, the \emph{finite element Hessian} and constructed a discrete
variational problem.

We proved that the numerical solution of this discrete variational
problem converges to the extrema of the continuous problem and that
the finite element Hessian converges to the Hessian of the continuous
extrema.

We foresee that this framework will prove useful when studying other
(possibly more complicated) second order variational problems, such as
discrete curvature problems like the affine maximal surface equation,
which is the topic of ongoing research.



\end{document}